\documentclass[10pt]{amsart}

\usepackage{amsmath,latexsym, amssymb,amsthm, amscd,graphicx}

\newtheorem{theorem}{Theorem}[section]
\newtheorem{proposition}[theorem]{Proposition}
\newtheorem{lemma}[theorem]{Lemma}
\newtheorem{corollary}[theorem]{Corollary}

\theoremstyle{remark}
\newtheorem{remark}[theorem]{Remark}
\theoremstyle{definition}

\def\real{\mathbb{R}}
\def\cone{\mathbf{C}}

\def\P{\mathcal{P}}
\def\m{\mathbf{m}}
\def\supp{\mathrm{supp}}

\def\A{\mathcal{A}} 
\def\a{\mathbf{a}}
\def\b{\mathbf{b}}
\def\c{\mathbf{c}}
\def\t{\mathbf{t}}
\def\Jac{Jac}
\def\base{\tau}
\def\cor{\mathbf{Cor}}
\def\ite{I}
\def\bbeta{\beta}
\def\gh{g}

\title[Decay of correlations]{Decay of correlations in suspension semi-flows \\ of angle-multiplying maps}
\author{Masato Tsujii}
\address{Department of Mathematics, Kyushu university, 
Fukuoka, 812-8581 Japan}

\email{tsujii@math.kyushu-u.ac.jp}

\date{\today}

\begin{document}
\begin{abstract}
We consider suspension semi-flows of angle-multiplying maps on the circle. Under a $C^r$generic condition on the ceiling function, we show that there exists an anisotropic Sobolev space\cite{BT} contained in the $L^2$ space such that the Perron-Frobenius operator for the time-$t$-map
acts on it and that the essential spectral radius of that action is bounded by the square root of the inverse of the minimum expansion rate.
This leads to a precise description on decay of correlations and extends the result of M. Pollicott\cite{Po}.
\end{abstract}
\maketitle

\section{Introduction}\label{sec:intro}
Decay of correlations and related topics for hyperbolic dynamical systems have been studied for more than three decades since the works of Bowen\cite{B}, Ruelle\cite{R} and Sinai\cite{S}.  
For the cases of discrete dynamical systems such as iterations of expanding maps and Anosov diffeomorphisms,
we nowadays have fairly good understanding on the rate of decay of correlations. (See \cite{VB, BT,GL1, GL2} and the references therein.) On the contrary, for the cases of continuous dynamical systems such as Anosov flows, the corresponding argument is much subtler and our knowledge is less satisfactory at present. 
A simple reason for the subtleness in the cases of flows is  that the time-$t$-maps of hyperbolic flows are {\em not} hyperbolic (but partially hyperbolic) as there is no expansion or contraction in the flow direction. 
It is only recent that Dolgopyat\cite{Do} showed rapid decay of correlations for topologically mixing Anosov flows and that Liverani\cite{L} showed exponential decay of correlations for Anosov flows preserving contact structures. Still we do not know, for instance, whether  
we observe  exponential decay of correations and quasi-compactness of the semi-groups of Perron-Frobenius operators for mixing (or generic) Anosov flows.

The aim of this paper is to study decay of correlations not for Anosov flows but for a class of expanding semi-flows,  suspension semi-flows of angle-multiplying maps on the circle, which we would like to view as a simplified model of the Anosov flow. 
We  consider Perron-Frobenius operators for the time-$t$-maps of such semi-flows and let them act on the anisotropic Sobolev spaces introduced in \cite{BT}. Our main result is that, under a $C^r$generic condition on the ceiling function,  the essential spectral radius of the action is bounded by the square root of the inverse of the minimum expansion rate.
This leads to a precise description on decay of correlations, which resembles the results known for hyperbolic discrete  dynamical systems\cite{BT,GL1}, and extends the earlier result of M.~Pollicott\cite{Po} on exponential decay.   

Actually a prototype of the argument in this paper has been appeared in \cite{AGT}, where a class of volume-expanding hyperbolic endomorphisms, called {\em fat solenoidal attractors}, were studied. 
In this paper, we will apply essentially the same idea to analyze the time-$t$-maps of the class of expanding semi-flows mentioned above. 
We emphasize that our intension is to display the idea in a simple setting 
 and present this paper as a study for the cases of hyperbolic flows\footnote{It will be possible to apply the idea presented in this paper to Anosov flows with smooth stable foliations. But we suspect that we  need to overcome essential  difficulties to treat more general Anosov flows.}. For this reason, we will confine our argument to a rather restrictive setting, though it is not very difficult to extend it to more general classes of expanding semi-flows and transfer operators. (See Remark \ref{re:generalization}.)

\begin{figure}[htbp]

\begin{center}
\includegraphics{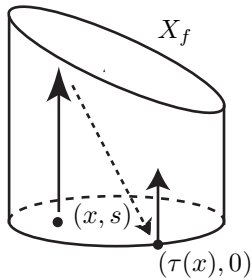}
\begin{picture}(0,0)
\put(-62,10){$(x,s)$}
\put(-30,-8){$(\tau(x),0)$}
\put(-30, 80){$X_f$}
\end{picture}
\caption{The semi-flow $\mathbf{T}_f$}
\label{fig:semiflow}
\end{center}
\end{figure}

We fix integers $\ell\ge 2$ and $r\ge 3$. Let $\base:S^1\to S^1$ be the angle-multiplying map on the circle $S^1=\mathbb{R}/\mathbb{Z}$ defined by $\base(x)=\ell x$. Let $C^r_+(S^1)$ be the space of positive-valued $C^r$functions on $S^1$. For each $f\in C^r_+(S^1)$, 
we consider the subset
\[
X_f=\{(x,s)\in S^1\times \real\mid 0\le s<f(x)\} 
\]
of the cylinder $S^1\times \real$.
The suspension semi-flow  
$\mathbf{T}_f=\{T^t_f:X_f\to X_f\}_{t\ge 0}$ of $\base$ is the semi-flow on $X_f$ in which each point on $X_f$ moves right upward (or the $s$-direction) with the unit speed and, at the instant it reaches to the upper boundary of $X_f$, it jumps down to the lower boundary with the $x$-coordinate transfered by $\tau$. (See Figure \ref{fig:semiflow}.)
The precise expression for its time-$t$-map is
\[
T^t_f(x,s)=\big(\base^{n(x,s+t;f)}(x), s+t-f^{(n(x,s+t;f))}(x)\big)
\]
where $f^{(n)}(x)=\sum_{i=0}^{n-1}f(\base^i(x))$ and $
n(x,t ;f)=\max\{n\ge 0 \mid  f^{(n)}(x)\le t\}$.

Let $m=m_f$ be the normalization of the restriction of the standard Lebesgue measure on $S^1\times \real$ to $X_f$. This is an ergodic invariant probability measure for $\mathbf{T}_f$. 
For a point $z=(x,s)\in X_f$ and $t\ge 0$, $E(z,t;f):=\ell^{n(x,s+t;f)}$ is 
the expansion rate along the orbit of $z$ up to time~$t$. 
The minimum expansion rate of $\mathbf{T}_f$ is naturally defined by
\[
\lambda_{\min}(\mathbf{T}_f)=\lim_{t\to \infty}\left(\min_{z\in X_f}E(z,t;f)\right)^{1/t}.
\]
For functions $\psi$ and $\varphi$ in $L^2(X_f)=L^2(X_f,m_f)$, we consider the correlation
\[
\cor_t(\psi,\varphi)=\int \psi\cdot \varphi\circ T^t_f 
 \;dm_f -\left(\int \varphi dm_f\right)
\left(\int \psi dm_f\right)\quad\mbox{for $t\ge 0$.}
\]  
Suppose that $\int \psi dm_f=0$ for simplicity.
If the semi-flow $\mathbf{T}_f$ is mixing,  we have  that $
\lim_{t\to \infty}\cor_t(\psi,\varphi)=0$. 
The question is the rate of convergence in this limit.  
In~\cite{Po}, M. Pollicott showed, under a mild condition on~$f$, that the rate is exponential: $
\left|\cor_t(\psi,\varphi)\right|<\mathrm{const.} \exp(-\epsilon t)$
for some $\epsilon>0$. (See also \cite{BV} for a generalization.) Under a $C^r$generic condition on~$f$, our results give a more precise description on asymptotic behavior of the correlation as $t\to \infty$ : For any real number $\mu>(\lambda_{\min}(\mathbf{T}_f))^{-1/2}$, there exists finitely many complex numbers\footnote{The complex numbers $\lambda_i$ may not be distinct from each other.} $\lambda_i\in \mathbb{C}$ with $\mu \le |\lambda_i|<1$ and integers $m_i\ge 0$ for $1\le i\le k$, such that
\begin{equation}\label{eq:exp}
\left|\;\cor_t(\psi,\varphi)-\sum_{i=1}^{k} H_i(\psi,\varphi)\cdot t^{m_i}\lambda_{i}^t \;\right|\le H_0(\psi,\varphi) \mu^t\quad \mbox{for $t\ge  0$}
\end{equation}
for any $\varphi\in L^2(X_f)$ and any $C^1$ function $\psi$ supported on the interior 
of $X_f$, where $H_i(\psi,\varphi)$ are coefficients that depend on $\psi$ and~$\varphi$ (and $\mu$).
As we will see later, the complex numbers $\lambda_i$ above are peripheral eigenvalues of the Perron-Frobenius operator for the time-$1$-map and each integer $m_i$ is bounded by the geometric multiplicity of the eigenvalue $\lambda_i$. 

In order to state the main results, we introduce some more notation. 
The differential $(DT^t_f)_z$ of $T^t_f$ at $z\in X_f$ is defined in the usual way if both $z$ and $T^t(z)$ belong to the interior of $X_f$ and, otherwise, is defined by 
\[
(DT^t_f)_z=\lim_{\epsilon\to +0} (DT^t_f)_{z+(0,\epsilon)}.
\] 
Fix\footnote{It is certainly better to choose $\gamma_0$ close to $1$, though it is not necessary.} a real number $\ell^{-1}<\gamma_0<1$.
Put 
\[
\theta_f=\max_{x\in S^1} |f'(x)|/(\gamma_0\ell-1)
\]
and 
\[
\cone_f=\cone(\theta_f):=\{(x,y)\in \real^2\mid |y|\le \theta_f |x|\}.
\] 
By the definition of $\theta_f$, we see that the cone  $\cone_f$ is strictly invariant for  $\mathbf{T}_f$ in the sense that
\[
(DT^t_f)_z(\cone_f)\subset \cone(\gamma_0\theta_f)\subset \cone_f
\quad\mbox{ for all $z=(x,s)\in X_f$ and $t\ge f(x)-s$.}
\]

Note that, for large $t$, the inverse image $(T_f^{t})^{-1}(z)$ of a point $z\in X_f$ consists of many points and thus there are many narrow cones $(DT^t_f)_{\zeta}(\cone_f)$ for $\zeta \in (T_f^{t})^{-1}(z)$ in the tangent space at $z$. 
As a measure for transversality between such cones, we introduce the quantity 
\[
\m(f,t)=\max_{z\in X_f}\max_{w\in (T^{t}_f)^{-1}(z)}
\sum_{\zeta:\zeta \not\pitchfork w}\frac{1}{E(\zeta,t;f)}
\]
where $\sum_{\zeta:\zeta \not\pitchfork w}$ is the sum over the points
$\zeta\in (T^{t}_f)^{-1}(z)$ such that
\[
(DT^t_f)_{\zeta}(\cone_f)\cap (DT^t_f)_w(\cone_f)\neq \{0\}.
\]
Note that we always have $\m(f,t)\le 1$ because 
\begin{equation}\label{sum}
\sum_{\zeta\in (T^{t}_f)^{-1}(z)}\frac{1}{E(\zeta,t;f)}= 1
\qquad \mbox{for any $t\ge 0$ and $z\in X_f$.}
\end{equation}
Finally we define the exponent
\[
\m(f)=\limsup_{t\to \infty} \m(f,t)^{1/t}\le 1.
\]

The Perron-Frobenius operator $\P^t_f:L^1(X_f)\to L^1(X_f)$ for $t\ge 0$ is defined by
\[
\P^{t}_f(u)(z)=\sum_{w\in (T^{t}_f)^{-1}(z)}\frac{u(w)}{\det{(DT^t_f)_w}},
\]
so that we have $\cor_t(\psi,\varphi)=\int \P^t_f \psi \cdot  \varphi\; dm_f$ provided that $\int \psi dm_f=0$. 

Let $C^1(X_f)$ be the set of functions $\varphi$ on $X_f$ such that $\P^t_f(\varphi)$ is $C^1$ on the interior of $X_f$ for any $t\ge 0$
\footnote{This condition imposes a restriction on the behavior of the function in the neighborhood of the  boundary of $X_f$, in addition to that it should be $C^1$ on the interior of $X_f$.}. This contains all the functions that are supported and $C^1$ on the interior of $X_f$.

Now the main results are stated as follows:
\begin{theorem}\label{th:main1}
There exists a Hilbert space $W_*(X_f)$ such that 
\[
C^{1}(X_f)\subset W_*(X_f)\subset L^{2}(X_f)
\]
and  that the Perron-Frobenius operator $\P^{t}_f$ for large $t\ge 0$ is restricted to the bounded operator $\P^{t}_f:W_*(X_f)\to W_*(X_f)$ whose essential spectral radius is bounded by 
$\m(f)^{t/2}$. 
\end{theorem}
\begin{theorem}\label{th:main2}
For each $\rho>1$,
there exists an open and dense subset $\mathcal{R}$ in  $C^r_+(S^1)$ such that, for $f\in \mathcal{R}$, the corresponding semi-flow\/ $\mathbf{T}_f=\{T^t_f\}$ is weakly mixing and satisfies
$\m(f)\le \rho\cdot \lambda_{\min}^{-1}(\mathbf{T}_f)$.
\end{theorem}
From these theorems, we obtain the following corollary.
\begin{corollary}\label{cor1}
For a $C^r$generic $f\in C^r_+(S^1)$, the semi-flow $\mathbf{T}_f$ is weakly mixing and the essential spectral radius of the Perron-Frobenius operator $\P^t_f$ acting on $W_*(X_f)$  is bounded by $\lambda_{\min}(\mathbf{T}_f)^{-t/2}$ for any sufficiently large $t$. 
\end{corollary}
The estimate (\ref{eq:exp}) for $C^r$ generic $f$ is an immediate consequence of this corollary. 
\begin{proof}[Proof of (\ref{eq:exp})]
Take large $t>0$. By Corollary \ref{cor1}, we have the decomposition $W_*(X_f)=E\oplus V$ where $E$ is the sum of the generalized eigenspaces for $\P^t_f$ corresponding to the eigenvalues not smaller than $\mu^t$ in absolute value  and 
\[
V=\{ \varphi \in W_*(X_f) \mid \mu^{-s} \|\P^s_f\varphi\|_* \to 0 \mbox{ as $s\to \infty$}\}
\]
where $\|\cdot\|_*$ denotes the norm on $W_*(X_f)$. 
The finite dimensional subspace $E$ is invariant under $\P^s_f$ for $s\ge 0$ by commutativity of $\P^s_f$ and $\P^t_f$. Thus we have $\P^s_f|_E=\exp(s B)$ for a linear map $B:E\to E$. Now it is easy to get the formula (\ref{eq:exp}), expressing $\psi$ as the sum of an element of $V$ and generalized eigenvectors of $B$. 
\end{proof}


Artur Avila told the author the following interesting observation related to the results above. Here we quote it under his permission. (But the author is responsible for the statement and proof.) 
\begin{theorem}\label{th:tfae} 
$\mathbf{T}_f$ is weakly mixing if and only if $\m(f)<1$.
\end{theorem}
This theorem and Theorem \ref{th:main1} (or Dolgopyat's argument\cite{Po}) imply that, once $\mathbf{T}_f$ is weakly mixing, it is exponentially mixing and there are no intermediate rates in  correlation decay (for functions in $C^1(X_f)$ at least). We will give the proof in Appendix~\ref{sec:apd1}. 

\begin{remark}\label{re:generalization}
We can generalize the argument in this paper to more general class of expanding semi-flows without much difficulty. For instance, the main results remains true with obvious changes in the related definitions when we consider arbitrary $C^3$ expanding maps on the circle in the place of $\tau$. 
The proof of Theorem \ref{th:main1} can be translated almost literally to such cases using the standard estimates on the distortion of expanding maps, while we need a bit modification to translate the proof of Theorem \ref{th:main2}. 
\end{remark}
Acknowledgement: The author would like to thank Artur Avila, S\'ebastien Gou\"ezel, Michihiro Hirayama and the anonymous referee for valuable comments that were crucial in improving this paper.


\section{Proof of Theorem \ref{th:main1}}
In this section we consider the semi-flow $\mathbf{T}_f=\{T^t_f\}_{t\ge 0}$ for some fixed $f\in C_+^r(S^1)$.
For simplicity, we write $T^t$ and $\P^t$ for $T^t_f$ and $\P^t_f$ respectively.

\subsection{Local charts on $X_f$}
We set up a system of local charts
\footnote{The system of local charts does {\em not} give the structure of (branched) manifold on $X_f$.} on $X_f$.  
To begin with, we consider two small real numbers $\eta>0$ and $\delta>0$, and set 
\[
R=(-\eta,\eta)\times (\delta,2\delta)
\subset
Q=(-2\eta,2\eta)\times (0,3\delta).
\]
For each $a=(x_0,s_0)\in X_f$ such that $[x_0-2\eta, x_0+2\eta]\times \{s_0\}\subset X_f$, we consider two mappings
\[
\kappa_{a}:Q\to X_f,  \quad \kappa_{a}(x,s)=T^s(x_0+x,s_0)
\]
and
\[
\tilde{\kappa}_{a}:Q\to S^1\times \real,  \quad \tilde{\kappa}_{a}(x,s)=(x_0+x,s_0+s).
\]
Note that $\kappa_a$ and $\tilde{\kappa}_a$ coincide when the image of $\tilde{\kappa}_a$ does not meet the upper boundary of $X_f$. 
Let $\eta$ and $\delta$ be so small that $\kappa_{a}$ is injective on $Q$ whenever it is defined.
 
Next we take a finite subset $A$ of $X_f$ so that the mapping $\kappa_a$ for $a\in A$ are defined and that the images $\tilde{\kappa}_{a}(R)$, $a\in A$, cover the subset  
\[
\tilde{X}_f:=\{ (x,s)\in S^1\times \real\mid
\delta/3\le s\le f(x) +2\delta/3\}.
\]
We may and do suppose that the intersection multiplicity of  
$\{\tilde{\kappa}_{a}(R)\}_{a\in A}$ is bounded by an absolute constant (say $100$). 
(This is in fact possible if we let the ratio $\eta/\delta$ be small.) 
Clearly the images of $\kappa_a(R)$ for $a\in A$ cover $X_f$. 

Let $C^r(R)$ be the set of $C^r$functions supported on $R$. We take the family $\{h_a\}_{a\in A}$ of functions in $C^\infty(R)$ as follows. First take a $C^\infty$ function $\bbeta_0:\real\to [0,1]$ such that $\bbeta_0(s)=1$ if $s\le \delta/3$ and $\bbeta_0(s)=0$ if $s\ge 2\delta/3$. We define $\bbeta:S^1\times \real\to [0,1]$ by
\[
\bbeta(x,s)=
\begin{cases}
\bbeta_0(s-f(x)), &\mbox{ if $s\ge f(x)$;}\\
1, &\mbox{ if $\delta<s< f(x)$;}\\
1-\bbeta_0(s), &\mbox{ if $s \le \delta$.}\\
\end{cases}
\]
This is a $C^\infty$function supported on $\tilde{X}_f$. From the choice of the finite subset $A$,  we can take a family $\{\tilde{h}_a:S^1\times \real\to [0,1]\}_{a\in A}$ of $C^{\infty}$ fucntions so that the support of each $\tilde{h}_a$ is contained in $\tilde{\kappa}_a(R)$ and that we have $\sum_a \tilde{h}_a\equiv \bbeta$ on $S^1\times \real$. We then define the $C^\infty$ function $h_a:\real^2\to [0,1]$ for $a\in A$ by 
\[
h_a=
\begin{cases}
\tilde{h}_a\circ \tilde{\kappa}_a &\mbox{ on $R$;}\\
0&\mbox{ on $\real^2\setminus R$.} 
\end{cases}
\]
\subsection{Anisotropic Sobolev spaces.}\label{ss:aniso}
We recall the  anisotropic Sobolev space and the and related definitions introduced in \cite{BT}. 
For a cone $\cone\subset \real^2$, we define its dual by
\[
\cone^*=\{v\in \real^2\mid (v,u)=0 \mbox{ for some $u\in \cone\setminus\{0\}$}\}.
\]
For two cones $\cone, \cone'\subset \real^2$, we write $\cone\Subset \cone'$  if the closure of $\cone$ is contained in the interior of $\cone'$ {\em except for} the origin.

A {\em polarization} $\Theta$ is a combination $\Theta=(\cone_+,\cone_-, \varphi_+, \varphi_-)$ of closed cones $\cone_{\pm}$ in $\real^2$ and $C^{\infty}$functions $\varphi_{\pm}:S^1\to [0,1]$ on the unit circle $S^1\subset \real^2$ that satisfy  
$\cone_+\cap \cone_-=\{0\}$
and 
\begin{equation}\label{vp}
\varphi_+(\xi)=
\begin{cases}
1, &\mbox{if $\xi\in S^1\cap \cone_{+}$;}\\
0, &\mbox{if $\xi\in S^1\cap \cone_{-}$,}
\end{cases} \qquad 
\varphi_-(\xi)=1-\varphi_+(\xi).
\end{equation} 
For two polarizations $\Theta=(\cone_+,\cone_-, \varphi_+, \varphi_-)$ and $\Theta'=(\cone'_+,\cone'_-, \varphi'_+, \varphi'_-)$, we write
$\Theta<\Theta'$ if $\real^2\setminus \cone'_+\Subset \cone_-$. 

Fix a $C^\infty$ function $\chi:\real\to [0,1]$ satisfying
\[
\chi(s)=
\begin{cases}
1 & \mbox{ for $s\le 1$;}\\
0 & \mbox{ for $s\ge 2$.}
\end{cases}
\]
For a polarization $\Theta=(\cone_+,\cone_-, \varphi_+, \varphi_-)$, an integer $n\ge 0$ and  $\sigma\in\{+,-\}$, we define the $C^{\infty}$function $\psi_{\Theta, n,\sigma}:\real^2\to [0,1]$ by
\[
\psi_{\Theta, n,\sigma}(\xi)=
\begin{cases}
\varphi_{\sigma}(\xi/|\xi|)\cdot (\chi(2^{-n}|\xi|)-\chi(2^{-n+1}|\xi|)),&\mbox{if $n\ge 1$;}\\
\chi(|\xi|)/2, &\mbox{if $n=0$.}
\end{cases} 
\]
This family of functions, $\psi_{\Theta, n,\sigma}$ for $n\ge 0$ and $\sigma\in \{+,-\}$, is a $C^\infty$ partition of unity on $\real^2$.

For a function $u\in C^{r}(R)$, we define  
\[
u_{\Theta, n,\sigma}(x)=\psi_{\Theta, n,\sigma}(D) u(x):=(2\pi)^{-2}\int e^{i(x-y)\xi} \psi_{\Theta, n,\sigma}(\xi) u(y) dy d\xi
\]
where $\psi_{\Theta, n,\sigma}(D)$ is the pseudo-differential operator with symbol $a(x,\xi)=\psi_{\Theta, n,\sigma}(\xi)$. Note that the pseudo-differential operator $\psi_{\Theta, n,\sigma}(D)$ may be viewed  as the composition $\mathcal{F}^{-1}\circ \Psi_{\Theta, n,\sigma}\circ \mathcal{F}$ where $\mathcal{F}$ is the Fourier transform and $\Psi_{\Theta, n,\sigma}$ is the multiplication operator by $\psi_{\Theta, n,\sigma}$.

For a polarization $\Theta=(\cone_+,\cone_-, \varphi_+, \varphi_-)$ and a real number $p$, we define the semi-norms $\|\cdot \|_{\Theta,p}^{+}$ and $\|\cdot \|_{\Theta,p}^{-}$ on $C^{r}(R)$ by
\[
\|u\|_{\Theta,p}^\sigma=\left(\sum_{n\ge 0} 2^{2pn}\|u_{\Theta, n,\sigma} \|^2_{L^2}\right)^{1/2}
=\left(\sum_{n\ge 0} 2^{2pn}\|\psi_{\Theta, n,\sigma} \cdot \mathcal{F}u\|^2_{L^2}\right)^{1/2}.
\]
Then we define the anisotropic Sobolev norm $\|\cdot \|_{\Theta,p,q}$ for real numbers $p$  and~$q$ by
\[
\|u\|_{\Theta,p,q}=\left((\|u\|_{\Theta,p}^+)^2
+(\|u\|_{\Theta,q}^-)^2
\right)^{1/2}.
\]
Clearly this norm is associated to a scalar product.

Actually we will not use the anisotropic Sobolev norms for general $p$, $q$ and $\Theta$ but those for the following specific cases. Fix small $0<\epsilon<1/2$ and put
\[
\|\cdot\|_{\Theta}^+:=\|\cdot\|_{\Theta,1}^+,\quad
\|\cdot\|_{\Theta}^-:=\|\cdot\|_{\Theta,0}^-,\quad
\|\cdot\|_{\Theta}:=\|\cdot\|_{\Theta,1,0}
\]
and
\[
|\cdot|_{\Theta}^+:=\|\cdot\|_{\Theta,1-\epsilon}^+,\quad
|\cdot|_{\Theta}^-:=\|\cdot\|_{\Theta,-\epsilon}^-,\quad
|\cdot|_{\Theta}:=\|\cdot\|_{\Theta,1-\epsilon,-\epsilon}.
\] 
In view of Parseval's identity, we have 
\begin{equation}\label{eqn:six}
\|\cdot \|_{\Theta}\ge \|\cdot\|_{L^2}/\sqrt{6}
\end{equation}
 where $6$ is (a bound for) the intersection multiplicity of the supports of $\psi_{\Theta,n,\sigma}$. 
The anisotropic Sobolev spaces $W_*(R;\Theta)$ and $W_\dag(R;\Theta)$ are the completion of $C^\infty(R)$ with respect to the norms $\|\cdot\|_{\Theta}$ and $|\cdot |_{\Theta}$ respectively. By (\ref{eqn:six}), we see that the space $W_*(R;\Theta)$ is naturally embedded in $L^2(R)$. Further we recall the next lemma from \cite{BT}: 
Let $W^s(R)$ be the (usual) Sobolev space of order $s$, that is, 
\[
W^s(R)=
\{ u\in \mathcal{D}'(\real^2)\;\mid\; \supp(u)\subset \mathrm{closure}({R}),\; (1+|\xi|^2)^{s/2}\mathcal{F}u(\xi)\in L^2(\real^2)\}.
\]  
\begin{lemma}\label{lm:cpt}For any polarizations $\Theta'<\Theta$, we have\\
{\rm{(a)}}\;$
C^1(R)\subset W^{1}(R)\subset W_*(R;\Theta)\subset L^2(R)$ and $ 
W^{1-\epsilon}(R)\subset W_\dag(R;\Theta)\subset W^{-\epsilon}(R)$.\\
{\rm{(b)}}\; $W_*(R;\Theta)\subset W_*(R;\Theta')$ and $W_\dag(R;\Theta)\subset W_\dag(R;\Theta')$.\\
{\rm{(c)}}
The inclusion $W_*(R;\Theta)\subset W_\dag(R;\Theta)$ is compact.
\end{lemma}
\begin{proof}
We can check (a) and (b) easily by using Parseval's identity. 
For the claim (c), we refer \cite[Proposition 5.1]{BT}. 
\end{proof}

As for the polarization $\Theta$, we will consider three  polarizations 
\[
\check{\Theta}_0=\big(\check{\cone}_{0,\pm}, \check{\varphi}_{0,\pm}\big)
<
\Theta_0=\big(\cone_{0,\pm}, \varphi_{0,\pm}\big)<
\hat{\Theta}_0=\big(\hat{\cone}_{0,\pm}, \hat{\varphi}_{0,\pm}\big)
\]
such that
\[
(\cone(\gamma_0\theta_f))^*\Subset
\hat{\cone}_{0,-}
\Subset
\left(\real^2\setminus \check{\cone}_{0,+}
\right) 
\Subset 
(\cone(\theta_f))^*.
\]

The Hilbert space $W_*(X_f)$ in Theorem \ref{th:main1} is defined as follows. 
Consider the projection operator (regarding functions as densities) 
\[
\Pi:(L^2(R))^A\to L^{2}(X_f), \quad \Pi((u_a)_{a\in A}) =\sum_{a\in A}\pi_a(u_a)
\]
where $\pi_a:L^2(R)\to L^2(X_f)$ for $a\in A$ is defined by
\[
\pi_a(u)(z)=
\begin{cases}
u(w)/\det D\kappa_a(w),&\quad\mbox{if $z=\kappa_a(w)$ for some $w\in R$;}\\
0, & \quad\mbox{otherwise.}
\end{cases}
\]
We equip the product spaces $(W_*(R;\Theta_0))^A\subset (L^2(R))^A$ and $(W_\dagger(R;\Theta_0))^A$  with the norms
\[
\|\mathbf{u}\|:=
\left(\sum_{a\in A}\|u_a\|_{\Theta_0}^2
\right)^{1/2} \mbox{ and }\quad
|\mathbf{u}|:=
\left(\sum_{a\in A}|u_a|_{\Theta_0}^2
\right)^{1/2}
\quad\mbox{where $\mathbf{u}=(u_a)_{a\in A}$}
\]
respectively,
so that they are Hilbert spaces. Then we put 
\[
W_*(X_f)=\Pi((W_*(R;\Theta_0))^A)
\]
 and equip it with the norm $
\|u\|=\inf \{\|\mathbf{u}\|\mid \Pi(\mathbf{u})=u\}$. 
This space $W_*(X_f)$ is isomorphic to the orthgonal complement of the kernel of $\Pi$ in $(W_*(R;\Theta_0))^A$ and, hence, is a Hilbert space.


\subsection{Transfer operators on local charts}
We consider the time-$t$-map $T^t$ and the corresponding Perron-Frobenius operator $\P^t$ for $t$ large enough\footnote{It is enough to consider $t$ with $t\ge \max_{x\in S^1}f(x)+3\delta$}. 
One technical difficulty in treating the time-$t$-map $T^t$ directly is that it is {\em not} continuous at the boundary of $X_f$. To avoid this problem, 
we consider a system of transfer operators on the local charts $\{\kappa_a\}_{a\in A}$ as a lift of $\P^t$, in which we do not find any trace of the discontinuity of $T^t$. 

For $a,b\in A$, let $
Q(a,b,t)$ be the set of points $z\in Q$ such that, for some $(x,s)\in Q$,
\begin{itemize}
\item[(i)] $T^t\circ \kappa_a(z)= \kappa_b((x,s))$  and 
\item[(ii)] $T^{t-\eta}\circ \kappa_a(z)\in \kappa_b(Q)$ for $0\le\eta<s$. 
\end{itemize}
Notice that the second condition (ii) does not follow from (i) when $\tilde{\kappa}_b(Q)$ meets the upper boundary of $X_f$. In fact, the second condition (ii) is important to avoid discontinuity in the following argument. 

For $a,b \in A$ with $Q(a,b,t)\neq \emptyset$, we consider the $C^r$ mapping
\[
T_{ab}^t:Q(a,b,t) \to Q, \quad T_{ab}^t(z)=\kappa_b^{-1}\circ T^t\circ \kappa_a(z).
\]
This is nothing but the mapping $T^t$ viewed in the local charts $\kappa_a$ and $\kappa_b$. 
Notice however that we restrict the domain of definition to $Q(a,b,t)$. 

For $a,b \in A$ with $Q(a,b,t)\neq \emptyset$, we consider the transfer operator 
\[
\P_{ab}:L^2(R)\to L^2(R),\quad \P_{ab}^tu(z)=
\sum_{w\in (T_{ab}^t)^{-1}(z)} \frac{h_b(z) u(w)}{\det (DT^t_{ab})_w}
\]
where the sum is taken over $w\in Q(a,b,t)$ such that $T_{ab}^t(w)=z$.
And we define the system of transfer operators on local charts,  
\[
\mathbf{P}^t:L^2(R)^{A}\to L^2(R)^{A},
\]
by 
\[
\mathbf{P}^t(\mathbf{u})=\left(\sum_{a\in A}\P^t_{ab}( u_{a})\right)_{b\in A} \quad \mbox{for $\mathbf{u}=(u_a)_{a\in A}\in L^2(R)^A$.}
\]
It is then not difficult to check that the following diagram commutes:
\begin{equation}\label{cd1}
\begin{CD}
L^1(R)^{A} @>{\mathbf{P}^t}>>L^1(R)^{A}\\
@VV{\Pi}V @VV{\Pi}V\\
L^1(X_f) @>{\P^t}>>L^1(X_f)
\end{CD}
\end{equation}

In the following subsections, we will prove  
\begin{proposition}\label{pp:LY}
The operator\/ $\mathbf{P}^t$ is restricted to the bounded operator
\[
\mathbf{P}^t:W_*(R;\Theta_0)^A\to W_*(R;\Theta_0)^A.
\]
Also $\mathbf{P}^t$ extends to the bounded operator 
\[
\mathbf{P}^t:W_\dag(R;\Theta_0)^A\to W_\dag(R;\Theta_0)^A.
\]
Further we have the Lasota-Yorke type inequality
\begin{equation}\label{eqn:LYP}
\|\mathbf{P}^t(\mathbf{u})\|\le C_\sharp\cdot  \m(f,t)^{1/2}
\|\mathbf{u}\|
+
C|\mathbf{u}|\qquad
\mbox{for $\mathbf{u}\in W_*(R)^A$}
\end{equation}
where the constant $C_\sharp$ does not depend on $t$ while the constant $C$ may. 
\end{proposition}
We can deduce Theorem \ref{th:main1} from this proposition and Lemma \ref{lm:cpt}(c). 

\begin{proof}[Proof of Theorem \ref{th:main1}]
Since $W_*(R;\Theta_0)^A$ is compactly embedded in $W_\dag(R;\Theta_0)^A$ from 
Lemma \ref{lm:cpt}(c),  the Lasota-Yorke type inequality in  
Proposition \ref{pp:LY} implies that 
 the essential spectral radius of the operator $\mathbf{P}^t:W_*(R;\Theta_0)^A\to W_*(R;\Theta_0)^A$ is bounded by $C_\sharp\cdot  \m(f,t)^{1/2}$. (See \cite{He, IM}.)
By the definition of the space $W_*(X_f)$, the commutative diagram (\ref{cd1}) is restricted to 
\begin{equation}\label{cd2}
\begin{CD}
W_*(R;\Theta_0)^{A} @>{\mathbf{P}^t}>>W_*(R;\Theta_0)^{A}\\
@VV{\Pi}V @VV{\Pi}V\\
W_*(X_f) @>{\P^t}>>W_*(X_f).
\end{CD}
\end{equation}
Recall also that the space $W_*(X_f)$ is identified with the orthogonal complement of the kernel of $\Pi$ in $W_*(R)^A$. Through this identification,  $\P^t$ corresponds to the composition of $\mathbf{P}^t$ with the orthogonal projection along the kernel of $\Pi$. Thus the essential spectral radius  of $\P^t:W_*(X_f)\to W_{*}(X_f)$ is bounded by that of $\mathbf{P}^t:W_*(R;\Theta_0)^A\to W_{*}(R;\Theta_0)^A$ or  $C_\sharp \m(f,t)^{1/2}$. 
Since this holds for any $t$ large enough, the essential spectral radius of $\P^t:W_*(X_f)\to W_{*}(X_f)$ is   bounded by $\m(f)^{1/2}$. The inclusions $C^1(X_f)\subset W_*(X_f)\subset L^2(X_f)$ follows from Lemma \ref{lm:cpt}(a). 
\end{proof}


\subsection{Two lemmas on the anisotropic Sobolev norms}
In this section, we give two basic lemmas on the anisotropic Sobolev norms. 
These lemmas and their proofs are slight modification of those given in 
\cite{BT}. For convenience of the reader, we give the proofs in Appendix \ref{sec:apd2} and \ref{sec:apd3}. (See also Remark \ref{qr}.)
\begin{lemma}\label{lm:pu}
Let $g_i:\real^2\to [0,1]$, $1\le i\le I$, be a family of $C^r$functions such that $\sum_{i=1}^{I}g_i(x)\le 1$ for $x\in Q$ and that $\supp(g_i)\subset Q$. Let $\Theta$ and $\Theta'$ be polarizations such that $\Theta'<\Theta$. Then we have
\begin{equation}\label{eqn:pu1}
\left[
\sum_{1\le i\le I}\|g_i u\|^2_{\Theta'}
\right]^{1/2}
\le C_0 \|u\|_{\Theta}+ C |u|_{\Theta}\qquad 
\mbox{for $u\in W_*(R;\Theta)$,  }
\end{equation}
where $C_0$ is a constant that does not depend on $\{g_i\}$ while the constant $C$ may. Further, if $\sum_{i=1}^{I}g_i(x)\equiv 1$ for all $x\in R$ in addition, we also have
\begin{equation}\label{eqn:pu2}
\|u\|_{\Theta'}\le \nu \left[
\sum_{1\le i\le I}\|g_i u\|^2_{\Theta}
\right]^{1/2}
+ C \sum_{i=1}^{I}|g_i u|_{\Theta}\qquad 
\mbox{for $u\in W_*(R;\Theta)$,}
\end{equation}
where $\nu$ is the intersection multiplicity of the supports of the functions $g_i$, $1\le i\le I$. 
\end{lemma}
In the next lemma, we consider the following situation. 
For a $C^{r-1}$function $h:\real^2\to \real$ supported on
a closed subset $K\subset R$ 
and for a $C^r$diffeomorphism $S:U\to S(U)\subset\real^2$ defined on an open neighborhood $U$ of $K$,   
we consider a transfer operator $L:C^{r-1}(R)\to C^{r-1}(R)$ defined by 
\[
Lu(x)=
\begin{cases}
h(x)\cdot  u\circ S(x),&\quad\mbox{ if $x\in K$;}\\
0&\quad\mbox{otherwise.}
\end{cases}
\]
Assume that, for polarizations $\Theta=(\cone_+,\cone_-, \varphi_+, \varphi_-)$ and $\Theta'=(\cone'_+,\cone'_-, \varphi'_+, \varphi'_-)$, we have  
\[
(DS_{\zeta})^{tr}(\real^2\setminus \cone_+)\Subset
\cone'_-\qquad \mbox{ for all $\zeta\in K$,}
\]
where $(DS_{\zeta})^{tr}$ denotes the transpose of $DS_{\zeta}$. 
Put 
\[
\gamma(S)=\min_{\zeta\in K} |\det DS_{\zeta}|
\]
and 
\[
\Lambda(S,\Theta',K)=\sup\left\{ \left. \frac{\|(DS_{\zeta})^{tr}(v)\|}{\|v\|}\; \right|\; \zeta\in K, \; 
(DS_{\zeta})^{tr}(v)\notin \cone'_-\;
\right\}.
\]
\begin{lemma}\label{lm:LYS}
The operator $L$ extends boundedly to $L:W_*(R;\Theta)\to
W_*(R;\Theta')$ and to $L:W_\dag(R;\Theta)\to
W_\dag(R;\Theta')$. Further we have, for $u\in W_*(R;\Theta)$,
\begin{equation}\label{eqn:LYSm}
\|Lu\|_{\Theta'}^-\le \gamma(S)^{-1/2} \|h\|_{L^\infty} \|u\|_{\Theta}
\end{equation}
and 
\begin{equation}\label{eqn:LYSp}
\|Lu\|_{\Theta'}^+\le C_0 \gamma(S)^{-1/2}  \Lambda(S,\Theta',K) \|h\|_{L^\infty} \|u\|_{\Theta}+C|u|_{\Theta}
\end{equation}
where the constant $C_0$ does not depend on $S$, $h$, $\Theta$ nor $\Theta'$ while the constant $C$ may. In particular, we have, for $u\in W_*(R;\Theta)$, 
\begin{equation}\label{eqn:LYSpm}
\|Lu\|_{\Theta'}\le C_0 \gamma(S)^{-1/2}\max\{1,\Lambda(S,\Theta',K)\} \|h\|_{L^\infty} \|u\|_{\Theta}+C|u|_{\Theta}.
\end{equation}
\end{lemma}
\begin{remark}\label{qr}
The latter claim (\ref{eqn:pu2}) of Lemma~\ref{lm:pu} is a special case of \cite[Lemma~7.1]{BT}. 
Also, Lemma \ref{lm:LYS} and the former claim (\ref{eqn:pu1}) in Lemma~\ref{lm:pu} correspond to \cite[Proposition~7.2 and 6.1]{BT}. However, since we considered only anisotropic Sobolev norms $\|\cdot\|_{\Theta,p,q}$ with $q<0<p$ in those propositions in \cite{BT},  we need to modify the statements and proofs slightly. 
This is the reason why we give the proofs of Lemma \ref{lm:LYS} and the former claim (\ref{eqn:pu1}) of Lemma~\ref{lm:pu} in the appendices \ref{sec:apd2} and \ref{sec:apd3}.  
\end{remark}

\def\E{\mathcal{E}}

\subsection{The Lasota-Yorke type inequality in local charts}
To complete the proof of Proposition \ref{pp:LY}, we have only  to show the Lasota-Yorke type inequality (\ref{eqn:LYP}), because the other claims now follow from Lemma \ref{lm:LYS} with $S$ the branches of the inverse of $T^t_{ab}$. For the proof of  the Lasota-Yorke type inequality (\ref{eqn:LYP}), it is enough to show the following lemma for the components $\P_{ab}^t$ of $\mathbf{P}^t$. 
\begin{lemma}\label{lm:LLY}
For $a,b\in A$ and for sufficiently large $t$, we have
\[
\|\P_{ab}^t(u)\|_{\Theta_0}^2\le C_\sharp\cdot  \m(f,t)
\|u\|_{\Theta_0}^2
+
C|u|_{\Theta_0}^2\qquad
\mbox{for $u\in C^r(R)$}
\]
where the constant $C_\sharp$ does not depend on $t$ while the constant $C$ may. 
\end{lemma}
Below we prove Lemma \ref{lm:LLY}, To begin with, we set up some notation. Fix $a,b\in A$ and consider large $t>0$. 
Let $\{D(\omega), \omega\in \Omega\}$, be a finite family of  small closed disks  whose interiors cover the closure of $R$. We may and do assume that the intersection multiplicity of this cover is bounded by some absolute constant (say $4$). Note that we can take such family of disks with arbitrarily small diameters. Let $D(\omega,i)$, $1\le i\le I(\omega)$, be the connected components of the preimage $(T_{ab}^{t})^{-1}(D(\omega))$ that meet the closure of~$R$. Then $\det DT^t_{ab}$ takes constant value on each component $D(\omega,i)$, which is denoted by $e(\omega,i)$. From (\ref{sum}), it holds
\begin{equation}\label{eqn:sum}
\sum_{1\le i\le I(\omega)} e(\omega,i)^{-1}\le 1 \quad\mbox{for any $\omega\in \Omega$.}
\end{equation}
We write $i\pitchfork_{\omega} j$ for $1\le i,j\le I(\omega)$ if 
\begin{equation}\label{eqn:ts}
(DT^t_{ab})_z(\cone_f)\cap (DT^t_{ab})_w(\cone_f)=\{0\}
\end{equation}
for any $z\in D(\omega,i)$ and $w\in D(\omega,j)$. 
Note that (\ref{eqn:ts}) holds if and only if
\[
(((DT^t_{ab})_z)^{tr})^{-1}((\cone_f)^*)\cap 
(((DT^t_{ab})_w)^{tr})^{-1}((\cone_f)^*)=\{0\}.
\]
Letting the diameters of the disks $D(\omega)$ be small, we may assume
\[
\sum_{j\not\pitchfork_\omega i} e(\omega,j)^{-1}\le  \m(f,t)
\quad \mbox{for $1\le i\le I(\omega)$ and $\omega\in \Omega$}
\]
where $\sum_{j\not\pitchfork_\omega i}$ denotes the sum over $1\le j\le I(\omega)$ such that $j\not\pitchfork_\omega i$. 
Further, we can 
take polarizations $\Theta(\omega,i)=(\cone_{\omega,i,+},\cone_{\omega,i,-}, \varphi_{\omega,i,+}, \varphi_{\omega,i,-})$ for each $\omega\in \Omega$ and $1\le i\le I(\omega)$ such that
\[
(((DT^t_{ab})_z)^{tr})^{-1})(\real^2\setminus \check{\cone}_{0,+})\Subset \cone_{\omega,i,-}\Subset 
\left(\real^2\setminus \cone_{\omega,i,+}\right)\Subset \hat{\cone}_{0,-}\quad
\mbox{for any $z\in D(\omega,i)$}
\]
and that 
\[
\overline{\left(\real^2\setminus \cone_{\omega,i,+}\right)} \cap \overline{\left(\real^2\setminus \cone_{\omega,j,+}\right)}=\{0\}\quad
\mbox{if $i\pitchfork_\omega j$.}
\]
Take a family of $C^\infty$ functions $\gh_\omega:\real^2\to [0,1]$ for $\omega\in \Omega$ such that the support of each $\gh_\omega$ is contained in $D(\omega)$ and that $\sum_{\omega\in \Omega}\gh_\omega(z)\equiv 1$ for all $z\in R$. We then define the functions $\gh_{\omega,i}:\real^2\to [0,1]$ for $\omega \in \Omega$ and $1\le i\le I(\omega)$ by 
\[
\gh_{\omega,i}(z)=
\begin{cases}
\gh_{\omega}(T_{ab}^{t}(z)), &\mbox{ if $z\in D(\omega,i)$ and}\\
0, &\mbox{ otherwise.}
\end{cases}
\] 
 
Now we start the proof of Lemma \ref{lm:LLY}. 
In the following, we will write $C_\sharp$ for constants that do not depend on $t$ while write $C$ for constants that may depend on~$t$. (Notice that the value of the constants denoted by $C_\sharp$ and $C$ are different from place to place.)  We view the operator $\P_{ab}^t$ under consideration as the composition of the four operations
\begin{itemize}
\item[(i)] breaking a function $u\in C^r(R)$ into $u_{\omega,i}:=\gh_{\omega,i} u$, $\omega\in \Omega$, $1\le i\le I(\omega)$, 
\item[(ii)] transforming  each $u_{\omega,i}$ to $v_{\omega,i}:=\P_{ab}^t(u_{\omega,i})$,  
\item[(iii)] summing up $v_{\omega,i}$ for $1\le i\le I(\omega)$ to get $v_{\omega}:=\sum_{1\le i\le I(\omega)} v_{\omega,i}=\gh_\omega \P^t_{ab}u$, 
\item[(iv)] summing up $v_{\omega}$ for $\omega\in \Omega$ to get $\P_{ab}^t u=\sum_{\omega}v_\omega$. 
\end{itemize}

For the operation (i), the former claim (\ref{eqn:pu1}) of Lemma \ref{lm:pu} gives the estimate
\[
\sum_{\omega\in \Omega}\sum_{i=1}^{I(\omega)}
\|u_{\omega,i}\|_{\check{\Theta}_0}^2\le (C_\sharp \|u\|_{\Theta_0} +C|u|_{\Theta_0})^2
\le C_\sharp \|u\|_{\Theta_0}^2 +C|u|^2_{\Theta_0}.
\]
For the operation (iv), the latter claim (\ref{eqn:pu2}) of Lemma \ref{lm:pu} gives the estimate 
\[
\|\P_{ab}^t u\|_{\Theta_0}^2 
=
\left\|\sum_{\omega\in \Omega} \gh_\omega \P_{ab}^t u
\right\|_{\Theta_0}^2=
\left\|\sum_{\omega\in \Omega} v_\omega\right\|_{\Theta_0}^2
\le 
C_\sharp \sum_{\omega\in \Omega}\|v_{\omega}\|_{\hat{\Theta}_0}^2+
C\sum_{\omega\in \Omega}|v_{\omega}|_{\hat{\Theta}_0}^2.
\]

Below we consider the operations (ii) and (iii). 
Letting $S$ be a branch of $(T_{ab}^{t})^{-1}$ and  $K=\supp(\gh_\omega)$ in Lemma \ref{lm:LYS}, we obtain the estimates
\begin{align}
&|v_{\omega,i}|_{\Theta(\omega,i)}\le 
C |u_{\omega,i}|_{\check{\Theta}_0}
\label{eqn:LY12}\qquad\mbox{and}\\
&\|v_{\omega,i}\|_{\hat{\Theta}_0}\le 
C_\sharp\cdot e(\omega,i)^{-1/2}\|u_{\omega,i}\|_{\check{\Theta}_0}
+ C|u_{\omega,i}|_{\check{\Theta}_0},
\label{eqn:LY11}\\
&\|v_{\omega,i}\|_{\hat{\Theta}_0}^+\le 
C_\sharp\cdot e(\omega,i)^{-3/2}\|u_{\omega,i}\|_{\check{\Theta}_0}+C|u_{\omega,i}|_{\check{\Theta}_0}.
\label{eqn:LY2}
\end{align}
The next lemma is the core of our argument, in which we make use of the transversality condition (\ref{eqn:ts})  essentially.  
\begin{lemma} \label{lm:trans}
If $i\pitchfork_\omega j$, we have
\begin{align*}
&\sum_{n\ge 0}\left|\left(\;\psi_{\hat\Theta_0, n,-}(D)v_{\omega,i}\;,\;\psi_{\hat\Theta_0, n,-}(D)v_{\omega,j}\;\right)_{L^2}\right|\le  C |v_{\omega,i}|_{\Theta(\omega,i)} |v_{\omega,j}|_{\Theta(\omega,j)}.
\end{align*}
\end{lemma} 
\begin{proof}
Put $w_{i,n}=\psi_{\hat\Theta_0, n,-}(D)v_{\omega,i}$, 
$w'_{i,n}= \psi_{\Theta(\omega,i),n,-}(D)v_{\omega,i}$,   
$w''_{i,n}=w_{i,n}-w'_{i,n}$. 
For $n>0$, we have, from the assumption $i\pitchfork_\omega j$, that
$
(w'_{i,n},  w'_{j,n})_{L^2}=0
$, that is,  
\[
(w_{i,n},w_{j,n})_{L^2}=
(w''_{i,n}, w'_{j,n})_{L^2}
+(w'_{i,n}, w''_{j,n})_{L^2}
+(w''_{i,n}, w''_{j,n})_{L^2}.
\]
Since we have that $\|w'_{i,n}\|_{L^2}\le  2^{\epsilon n}|v_{\omega,i}|_{\Theta(\omega,i)}$ and that $\|w''_{i,n}\|_{L^2}\le  2^{-(1-\epsilon)n}|v_{\omega,i}|_{\Theta(\omega,i)}$ by the definition of the norm $|\cdot|_{\Theta(\omega,i)}$ and since $0<\epsilon<1/2$, we can get the lemma using Schwarz inequality.
\end{proof}

It follows from Lemma \ref{lm:trans}  that 
\[
\left(\left\|v_\omega
\right\|_{\hat\Theta_0}^{-}\right)^2
\le \sum_{i}\sum_{j\not\pitchfork_{\omega} i} \|v_{\omega,i}\|_{\hat\Theta_0}
\|v_{\omega,j}\|_{\hat\Theta_0}
+
C\sum_{i}|v_{\omega,i}|_{\Theta(\omega,i)}^2
\]
for some constant $C>0$ that may depend on $I(\omega)$, where $\sum_{j\not\pitchfork_{\omega} i}$ denotes the sum over $1\le j\le I(\omega)$ such that $j\not\pitchfork_{\omega} i$.
Applying the inequality (\ref{eqn:LY11}) and (\ref{eqn:LY12}) in the right hand side along, we obtain
\begin{align*}
\left(\left\|v_\omega
\right\|_{\hat\Theta_0}^{-}\right)^2
&\le C_{\sharp}\sum_{i}\sum_{j\not\pitchfork_{\omega} i}
\frac{e(\omega,j)^{-1} \|u_{\omega,i}\|_{\check\Theta_0}^2
+
e(\omega,i)^{-1}
\|u_{\omega,j}\|_{\check\Theta_0}^2}{2}
+
C\sum_{i}|u_{\omega,i}|_{\check\Theta_0}^2\\
&\le C_{\sharp}\m(f,t)\sum_{i}\|u_{\omega,i}\|_{\check\Theta_0}^2
+
C\sum_{i}|u_{\omega,i}|_{\check\Theta_0}^2.
\end{align*}
On the other hand we obtain, from (\ref{eqn:LY2}), that
\begin{align*}
\left(\left\|v_\omega
\right\|_{\hat\Theta_0}^{+}\right)^2
&\le C_\sharp \left(\sum_{i} e(\omega,i)^{-3/2} \|u_{\omega,i}\|_{\check\Theta_0}\right)^2
+
C\left(\sum_{i}|u_{\omega,i}|_{\check\Theta_0}\right)^2\\
&\le C_{\sharp}\m(f,t)\sum_{i} \|u_{\omega,i}\|_{\check\Theta_0}^2
+
C\sum_{i}|u_{\omega,i}|_{\check\Theta_0}^2
\end{align*}
where we used Schwarz inequality, (\ref{eqn:sum}) and the simple fact $e(\omega,i)^{-1}\le \m(f,t)$ in the second inequality. 
Therefore we obtain, for the operations (ii) and (iii),
\[
\left(\left\|v_\omega
\right\|_{\hat\Theta_0}\right)^2
\le C_{\sharp}\m(f,t)\sum_{i} \|u_{\omega,i}\|_{\check\Theta_0}^2
+
C\sum_{i}|u_{\omega,i}|_{\check\Theta_0}^2.
\]

Letting $S$ be the identity map in Lemma \ref{lm:LYS}, we see that $|u_{\omega,i}|_{\check\Theta_0}\le C|u|_{\Theta_0}$ and that $|v_{\omega,i}|_{\hat\Theta_0}\le C|v_{\omega,i}|_{\hat\Theta(\omega,i)}$.  These and (\ref{eqn:LY12}) give 
\[
|v_{\omega}|_{\hat\Theta_0}\le C\sum_{1\le i\le I(\omega)} |v_{\omega,i}|_{\Theta(\omega,i)}
\le C \sum_{1\le i\le I(\omega)} |u_{\omega,i}|_{\check\Theta_0}\le C|u|_{\Theta_0}.
\]
We can now conclude the Lasota-Yorke type inequality in Lemma \ref{lm:LLY} from the estimates on the operations (i)--(iv) given above. 


\section{Proof of Theorem \ref{th:main2}}
The proof of Theorem \ref{th:main2} below is a modification of the argument in~\cite{T}. 

\subsection{Notations}
We recall some notation from \cite{T}. 
Let $\A=\{1,2,\dots,\ell\}$ and let $\A^n$ be the space of words of length $n$ on $\A$.
For a word $\a=(a_i)_{i=1}^{n}\in \A^n$ and an integer $0\le p\le n$, let $[\a]_p=(a_i)_{i=1}^{p}$.
For $0\le p\le n$, we define the equivalence relation $\sim_p$ on $\A^n$ so that $\a\sim_p\b$ if and only if $[\a]_p=[\b]_p$.

Let $\P$ be the partition of $S^1$ into the intervals $\P(k)=[(k-1)/\ell, k/\ell)$ for $k\in \A$. Then $\P=\bigvee_{i=0}^{n-1}\base^{-i}(\P)$ is the partition into the intervals
\[
\P(\a)=\bigcap_{i=0}^{n-1}\base^{-i}(\P(a_{n-i})),\quad \a=(a_i)_{i=1}^{n}\in \A^n.
\] 
Let  $x_{\a}$ the left end point of the interval $\P(\a)$.
\begin{remark}
Notice that $\a$ is the {\em inverse} of the itinerary of the points in $\P(\a)$. 
\end{remark}

For a point $x\in S^1$ and $\a=(a_i)_{i=1}^{n}\in \A^n$, we denote by $\a(x)$ the unique point $y\in \P(\a)$ such that $\base^{n}(y)=x$.  For a $C^r$function $f\in C^r(S^1)$, $x\in S^1$ and $\b\in \A^n$, we put
\[
s(x,\b;f):=f^{(n)}(\b(x))=\sum_{i=1}^{n} f([\b]_i(x)).
\]
Then we have 
\[
\frac{ds}{dx}(x,\b;f)=\sum_{i=1}^{n}\ell^{-i} \frac{df}{dx}([\b]_i(x)).
\]

We will identify the unit circle $S^1$ with the lower boundary $S^1\times\{0\}$ of  $ X_f$. If $f^{(n)}(\b(x))\le t<f^{(n+1)}(\b(x))$, the image of the horizontal tangent vector $(1,0)$ at $\b(x)\in S^1\times\{0\}$ by the mapping $T^{t}_f$ has slope  $\frac{d}{dx}s(x,\b;f)$ and hence
\[
(DT^t_f)_{\b(x)}(\cone_f)=
\left\{
(\xi,\eta)\in \real^2\;\left|\;\;\;  
\left|\eta-
\frac{ds}{dx}(x,\b;f)\cdot \xi \right|
\le \ell^{-n}\theta_f |\xi|
\right.
\right\}.
\]

For $K>1$, let $C^r_+(S^1;K)$ be the set $f\in C^r_+(S^1)$ such that $K^{-1}< f(x)< K$ for $x\in S^1$ and $\|f\|_{C^r}< K$. 
By virtue of Theorem \ref{th:tfae}, it is sufficient for the proof of Theorem \ref{th:main2} to show that, for each $\rho>1$ and $K>0$,  
the condition 
\begin{equation}\label{eqn:claimrho}
\m(f)\le \rho\cdot \lambda_{\min}(\mathbf{T}_f)^{-1}
\end{equation}
holds for functions $f$ in an open and dense subset of $C^r_+(S^1;K)$. 
We will prove this claim in the following. We henceforth fix arbitrary $\rho>1$ and  $K>0$. Note that we have, for  $f \in C^r_+(S^1;K)$, that
\[
\theta_f\le \theta_K:= K/(\gamma_0\ell-1),
\qquad \ell^{1/K}\le \lambda_{\min}(\mathbf{T}_f)\le \ell^K
\]
and also that
\begin{equation}\label{eqn:tk}
\left|\frac{d^2s}{d^2x}(x,\b;f)\right|
=\left|\sum_{i=1}^{n}\ell^{-2i} \frac{d^2f}{dx^2}([\b]_i(x))\right|
\le \frac{\ell^{-2}K}{1-\ell^{-2}}\le \theta_K
\end{equation}
for $x\in S^1$ and $\b\in \A^n$.


\subsection{Some consequences of the condition $\m(f)> \rho\cdot \lambda_{\min}(\mathbf{T}_f)^{-1}$} 
In this subsection, we see what kind of singular situation occurs if the condition  (\ref{eqn:claimrho}) does not hold. 
Fix $1<\gamma<\ell$ such that $\gamma^K<\rho$.   
\begin{proposition}\label{pp:a}
If $\m(f)> \rho\cdot \lambda_{\min}(\mathbf{T}_f)^{-1}$ for $f\in C^r_+(S^1;K)$, then,  for any  $n\ge 1$, there exist $\c\in \A^{n}$ and $B\subset \A^{n}$ with $\# B\ge  \gamma^{n}$ such that 
\[
\left|\frac{ds}{dx}(x_{\c},\a;f)-\frac{ds}{dx}(x_{\c},\b;f)\right|\le 8\theta_K\cdot  \ell^{-n}
\quad\mbox{for all $\a,\b\in B$.}
\]
\end{proposition}
\begin{proof} 
Take $\gamma<\bar\gamma<1$ so close to $\gamma$ that $\bar\gamma^K<\rho$.  Then take $1<\lambda<\lambda_{\min}(\mathbf{T}_f)$ so close to $\lambda_{\min}(\mathbf{T}_f)$  that $\bar\gamma^K\lambda_{\min}(\mathbf{T}_f)<\rho\lambda$. From the assumption, we can take an arbitrarily large $t\ge 0$, $z\in X_f$, $w\in T^{-t}(z)$ and $Z\subset T^{-t}(z)$ such that 
\begin{equation}\label{eqn:tang}
(DT^t)_\zeta(\cone_f)\cap (DT^t)_w(\cone_f)\neq \{0\}
\quad\mbox{ for $\zeta\in Z$}
\end{equation}
and that 
\begin{equation}\label{eqn:sz}
\sum_{\zeta\in Z} \frac{1}{E(\zeta,t;f)}\ge \rho^t\cdot \lambda_{\min}(\mathbf{T}_f)^{-t}.
\end{equation}
We may and do assume in addition that $z\in S^{1}\times\{0\}$. Put 
\[
m=\min\left\{n(x,s+t;f)\;\left|\; (x,s)\in Z\right.\right\}.
\]
 (Recall the definition of $n(x,t;f)$ in Section \ref{sec:intro}.) 
Then we have 
\begin{equation}\label{qnn:Kb}
K^{-1}t<m<Kt\quad \mbox{ and }\quad \ell^{m}> \lambda^{t}
\end{equation}
provided that $t$ is sufficiently large. 
(To see the second inequality, note that we have $\ell^m\sim \lambda_{\min}(\mathbf{T}_f)^{t}$ for large $t$.)

For each $\zeta=(x,s)\in Z$, let  $\ite(\zeta)\in \A^{n(x,s+t;f)}$ be the word such that $\zeta\in \P(\ite(\zeta))$.
Put $A=\{[I(\zeta)]_{m}\mid \zeta\in Z\}\subset \A^{m}$. 
Then, applying the formula (\ref{sum}), we can see that
\[
\# A\cdot \ell^{-m}\ge
\sum_{\zeta\in Z} \frac{1}{E(\zeta,t;f)}.
\]
Combining this inequality with (\ref{eqn:sz}), we see
\begin{equation}\label{eqn:cardA}
\# A\ge \ell^{m} \rho^{t}\lambda_{\min}(\mathbf{T}_f)^{-t}>(\lambda\rho \cdot \lambda_{\min}(\mathbf{T}_f)^{-1})^t 
> \bar\gamma^{Kt}>
\bar\gamma^m.
\end{equation}
Note that the condition (\ref{eqn:tang}) implies 
\begin{equation}\label{eqn:difds}
\left|\frac{ds}{dx}(z, \b;f)-\frac{ds}{dx}(z, \b';f)\right|\le 4\theta_K\cdot \ell^{-m}
\quad\mbox{ for all $\b,\b'\in A$.}
\end{equation}

Now we start to 
consider an arbitrary integer $n\ge 1$ in the claim. Taking large $t$ in the beginning, we may and do assume that the integer $m$ above is much larger than~$n$. 
For each $0\le k\le m$, 
we split $A$ into equivalence classes with respect to  $\sim_k$ and 
let $A_k\subset A$ be one of those equivalence classes with 
maximum cardinality. 
Then the cardinality  $q(k)$ of $A_k$ is a decreasing sequence with respect to $k$. 
Further we have $q(0)\ge\bar\gamma^m$ by (\ref{eqn:cardA}) and also $q(m)=1$ obviously. Therefore we can find an integer $0\le m'\le m-n$ such that $q(m'+n)\le \gamma^{-n} q(m')$, provided that we took sufficiently large $t$ in the beginning. 
We fix such integer $m'$ and let  $A'\subset \A^{m-m'}$ be the set of words that are obtained by removing the first common $m'$ letters (say $\c$) from the words in $A_{m'}$, and put $x=\c(z)$. Then it follows from (\ref{eqn:difds}) that 
\[
\left|\frac{ds}{dx}(x, \b;f)-\frac{ds}{dx}(x, \b';f)\right|\le 4\theta_K\cdot\ell^{-(m-m')}\le 4\theta_K\cdot\ell^{-n}
\quad\mbox{ for all $\b,\b'\in A'$.}
\]

Put $B=\{[\a]_{n}\mid \a\in A'\}$. From the condition $q(m'+n)\le \gamma^{-n} q(m')$ in the choice of $m'$, we have that $\# B\ge \gamma^n$. Also, since 
\[
\left|\frac{ds}{dx}(x,[\a]_n ;f)-\frac{ds}{dx}(x, \a;f)\right|\le \theta_K\cdot\ell^{-n}
\quad\mbox{ for $\a \in \A^{m-m'}$, }
\]
we see that
\[
\left|\frac{ds}{dx}(x, \b;f)-\frac{ds}{dx}(x, \b';f)\right|\le 6\theta_K\cdot\ell^{-n}
\quad\mbox{ for all $\b,\b'\in B$.}
\]

For $\b,\c\in \A^n$ and $f\in C_+^r(S^1;K)$, the variation of  $\frac{d}{dx}s(\cdot,\b;f)$ on the interval $\P(\c)$ is bounded by 
$\theta_K \ell^{-n}$, in view of (\ref{eqn:tk}). Therefore, translating the point $x$ to the point $x_\c$, we obtain the conclusion of the proposition. 

\end{proof}

\def\d{\mathbf{d}}

To state the next proposition, we set up some constants. First 
take real numbers $\alpha$ and $\beta$ such that $1<\beta<\alpha<\gamma$  and then take positive integers $p$ and $\nu$ such that
\[
\beta^{-p}\ell^2<1\quad
\mbox{and}\quad
(\nu+1)(p+1)\alpha^{-\nu}<1.
\]
We put 
\[
\delta=\frac{\log \gamma -\log \alpha}{\log \ell-\log \alpha} \;\in (0,1).
\]
Then we  choose an integer $N>\nu$ such that 
\[
\ell^\nu\alpha^{n}<\gamma^n\quad \mbox{for $n\ge N$}
\]  
and that 
\[
\ell^{-\nu}(\gamma/\beta)^{n'}(1-(\nu+1)(p+1) \alpha^{-\nu})\ge 1 \quad \mbox{for $n'\ge \delta N$.  }
\]
\begin{proposition}\label{pp:aa}
If $\m(f)> \rho\cdot \lambda_{\min}(\mathbf{T}_f)^{-1}$ for $f\in C^r_+(S^1;K)$, then, for any $n\ge N$,  there exist an integer $\delta n\le n'\le n$,  a word $\d\in \A^{n'}$ and mutually disjoint subsets $B_i\subset \A^{n'}$ for $1\le i\le (\nu+1)(p+1)$ such that 
\begin{itemize}
\item[\rm{(a)}]
$\displaystyle \left|\frac{ds}{dx}(x_{\d},\b;f)-\frac{ds}{dx}(x_{\d},\b';f)\right|\le 10 \theta_K\cdot \ell^{-n'}$ for all $\b,\b'\in \bigcup_{i=1}^{(\nu+1)(p+1)} B_i$,
\item[\rm{(b)}]$\# B_i\ge  \beta^{n'}$ for $1\le i\le  (\nu+1)(p+1)$, and 
\item[\rm{(c)}]
$[\a]_\nu=[\b]_\nu$ for $\a\in B_i$ and $\b\in B_j$  if and only if $i=j$.
\end{itemize} 
\end{proposition}
\begin{proof}
Let $n\ge N$ and let $B\subset \A^n$ be the subset in the conclusion of Proposition~\ref{pp:a}. 
For $0\le k\le [n/\nu]$, we split $B$ into equivalence classes with respect to $\sim_{k\nu}$ and 
let $B_k\subset B$ be one of those equivalence classes with maximum cardinality.  
Then the cardinality $q(k)$ of $B_k$ is decreasing with respect to $k$ and satisfies $q(0)=\# B \ge \gamma^n$ and $q([n/\nu])\le \ell^\nu< \gamma^n\alpha^{-[n/\nu]\nu}$, where the last inequality follows from the first condition in the choice of $N$.

Let $k_0$ be the smallest integer $1\le k\le [n/\nu]$ such that 
$q(k)< \gamma^n\alpha^{-k\nu}$. By this choice of $k_0$, we have 
\begin{equation}\label{eqn:qk0}
q(k_0)< \alpha^{-\nu} q(k_{0}-1)
\quad \mbox{ and } \quad q(k_{0}-1)\ge \gamma^{n}\alpha^{-(k_0-1)\nu}.
\end{equation}
 Put $n'=n-(k_0-1)\nu$. 
Since $q(k)\le \ell^{n-\nu\cdot k}$ obviously, we have 
\[
\ell^{n'}=\ell^{n-(k_0-1)\nu} \ge q(k_0-1)\ge \gamma^n\alpha^{n'-n}\quad \mbox{or} 
\quad n'\ge \frac{\log\gamma-\log \alpha}{\log \ell-\log \alpha}\cdot n=\delta n.
\]

Let $B'_i\subset B_{k_0 -1}$, $1\le i\le \ell^\nu$, be the equivalence classes in $B_{k_0 -1}$ with respect to the relation $\sim_{k_0}$, arranged in decreasing order of cardinality\footnote{Some of $B'_i$'s may be empty.}. 
Then we have a simple inequality
\begin{align*}
\min_{1\le i\le (\nu+1)(p+1)}\# B'_i
&\ge \frac{q(k_0-1)-(\nu+1)(p+1)q(k_{0})}{\ell^\nu}\\
&\ge \ell^{-\nu}\gamma^n\alpha^{-(k_0-1)\nu}(1-(\nu+1)(p+1) \alpha^{-\nu})\ge \beta^{n'}
\end{align*}
where the second inequality follows from (\ref{eqn:qk0}) and the last  from the second condition in the choice of $N$. 
Finally let $B_i\subset \A^{n'}$ for $1\le i\le (\nu+1)(p+1)$ be the subset of words that are obtained by 
removing the first common $(k_0-1)\nu$ letters (say $\c'$) from the words in $B'_i$. Then the conditions (b) and (c)  hold. 
From the condition on the subset $B$ in Proposition \ref{pp:a}, we have
\[
\left|\frac{ds}{dx}(x_{\c\c'},\b;f)-\frac{ds}{dx}(x_{\c\c'},\b';f)\right|\le 8 \theta_K\cdot \ell^{-n'}\quad\mbox{ for all $\b,\b'\in \bigcup_{i=0}^{(\nu+1)(p+1)} B_i$.}
\]
Take $\d\in \A^{n'}$ such that $x_{\c\c'}\in \mathcal{P}(\d)$.
($\c$ is that in Proposition \ref{pp:a}.) Then the condition (a) holds because the variations of the functions $\frac{d}{dx}s(\cdot,\a;f)$ for $\a\in \A^{n'}$ on $\P(\d)$ are bounded by 
$\theta_K \ell^{-n'}$, in view of (\ref{eqn:tk}).
\end{proof}

\subsection{Generic perturbations}
We are going to show that the consequences of the condition $\m(f)> \rho\cdot \lambda_{\min}(\mathbf{T}_f)^{-1}$ given is Proposition \ref{pp:aa} hold only for very small set of $f\in C^r_+(S^1;K)$. For this purpose, we next consider about perturbations of the function $f$. 

For $f\in C_+^r(S^1;K)$ and $\varphi_i\in C^{\infty}(S^1)$, $1\le i\le m$,  we consider the family
\begin{equation}\label{eqn:fam}
f_{\t}(x)=f(x)+\sum_{i=1}^{m}t_i\cdot \varphi_i(x)
\end{equation}
with parameter $\t=(t_i)_{i=1}^{m}\in \real^m$.
For a point $x\in S^1$ and a finite subset $\sigma=\{\b_i\}_{0\le i\le p}$ of $\A^n$, let  $G_{x,\sigma}:\real^m\to \real^{p}$ be the affine map defined by
\[
G_{x,\sigma}(\t)=
\left(
\frac{ds}{dx}(x,\b_i;f_{\t})-\frac{ds}{dx}(x,\b_0;f_{\t})
\right)_{i=1}^{p}.
\]
Note that $G_{x,\sigma}(\t)$ is independent of $f$ in (\ref{eqn:fam}).  
For an affine map $A:\real^m\to \real^p$, let $\Jac(A)$ be the Jacobian of $DA|_{\ker(DA)^{\perp}}$, the restriction of the linear part $DA$ to the orthogonal complement of its kernel when $A$ is surjective, and put $\Jac(A)=0$ otherwise. 
In other words, $\Jac(A)$ is 
the maximum among the Jacobians of the restrictions of $DA$ to $p$-dimensional subspaces in $\real^m$. 
The following is a slight variant of \cite[Proposition 16]{T}. 
\begin{proposition}\label{pp:b}
We can choose functions $\varphi_i\in C^{\infty}(S^1)$, $1\le i\le m$, 
such that, for any $x\in S^1$ and any subsets 
$A=\{\a_i\}_{1\le i\le (\nu+1)(p+1)}$ of $\A^\nu$, there exist a subset $A'=\{\a'_i\}_{0\le i\le p}$ of $A$ such that we have $\Jac(G_{x,\sigma})\ge 1$ whenever a subset $\sigma=\{\b_i\}_{0\le i\le p}$ of $\A^n$ with $n\ge \nu$ satisfies  $[\b_i]_{\nu}=\a'_i$ for $0\le i\le p$.
\end{proposition}
The proof of Proposition \ref{pp:b} is similar to that of  \cite[Proposition 16]{T}. For completeness, we give the proof in the last subsection.


\subsection{The end of the proof}
For $n\ge \nu$, $\c\in \A^n$ and $\sigma=(\b_i)_{i=0}^{p}\in (\A^n)^{p+1}$, let $Y(n,\c,\sigma)$ be the set of functions $f\in C^r_+(S^1;K)$ such that
\[
\left|\frac{ds}{dx}(x_{\c},\b_i;f)-\frac{ds}{dx}(x_{\c},\b_0;f)\right|\le 10\theta_K\cdot \ell^{-n}
\quad\mbox{for all $1\le i\le p$.}
\]
Note that $Y(n,\c,\sigma)$ is a closed subset in $C^r_+(S^1;K)$.

For $n\ge \nu$, let $Y(n)$ be the set of  functions $f\in C_+^r(S^1;K)$ that belongs to 
$Y(n,\c,\sigma)$ for more than $[\beta^{n(p+1)}]$ combinations of $(\c,\sigma)\subset \A^n\times  (\A^n)^{p+1}$ satisfying $\Jac(G_{x_{\c},\sigma})\ge 1$. Let $Y_*(n)=\bigcup_{n'=[\delta n]}^{n} Y(n')$. Then $Y(n)$ and $Y_*(n)$ are also closed subsets in $C^r_+(S^1;K)$. 
Proposition \ref{pp:aa} tells that, if $\m(f)> \rho\cdot \lambda_{\min}(\mathbf{T}_f)^{-1}$, then $f$ belongs to the closed subset $ \bigcap_{n\ge N}Y_*(n)$. To finish the proof of the theorem,  we show that the complement of $\bigcap_{n\ge N}Y_*(n)$ is dense in $C^r_+(S^1;K)$.

Take a function $f\in C^r_+(S^1;K)$ arbitrarily and consider the family (\ref{eqn:fam}) with $\varphi_i\in C^{\infty}(S^1)$, $1\le i\le m$, in Proposition \ref{pp:b}. 
Take $\epsilon>0$ so small that $f_\t\in C_+^r(S^1;K)$ for all $\t\in [-\epsilon,\epsilon]^m$. 
Let $X(n,\c,\sigma)$, $X(n)$ and $X_*(n)$ be the set of parameters $\t\in [-\epsilon,\epsilon]^m$ such that $f_\t\in Y(n,\c,\sigma)$, that $f_\t\in Y(n)$ and that $f_\t\in Y_*(n)$, respectively. 
From the definition of Jacobian in the last subsection, we have $
Leb(X(n,\c,\sigma))\le C \ell^{-np}$
for some constant $C>0$ that depends on $\theta_K$, $m$ and $\epsilon$.
Therefore, taking the number of combinations of $(\c,\sigma)$ into consideration, we get
\[
Leb(X(n))\le \frac
{C\ell^{-np}\times \ell^n\times \ell^{(p+1)n}}
{\beta^{(p+1)n}}<C(\beta^{-p}\ell^2)^n. 
\]
As we chose $p$ such that $\beta^{-p}\ell^2<1$, we have 
$Leb(\bigcap_{n\ge N} X_*(n))=0$ and hence the complement of $\bigcap_{n\ge N} Y_*(n)$ in $C^r_+(S^1;K)$ is dense. 

\begin{remark}
The proof above shows also that the condition $\m(f)\le \lambda_{\min}(\mathbf{T}_f)^{-1}$ holds for a prevalent subset of $f\in C^r_+(S^1)$ in measure-theoretical sense(\cite{HSY, T2}).
\end{remark}

\subsection{The proof of Proposition \ref{pp:b}}
To prove Proposition \ref{pp:b}, it is enough to show the following localized version of the claim.
\begin{proposition}\label{pp:b2}
For each $y\in S^1$, we can choose functions $\varphi_{y,i}\in C^{\infty}(S^1)$ for $1\le i\le \ell^\nu$ and a neighborhood $U_y$ of $y$ such that, for any point $x\in U_y$ and any subsets 
$A=\{\a_i\}_{1\le i\le (\nu+1)(p+1)}$ of $\A^\nu$, there exists a subset $A'=\{\a'_i\}_{0\le i\le p}$ of $A$ such that we have that $\Jac(G_{x,\sigma})\ge 1$ whenever a subset $\sigma=\{\b_i\}_{0\le i\le p}$ of $\A^n$ with $n\ge \nu$ satisfies $[\b_i]_{\nu}=\a'_i$ for $0\le i\le p$.
\end{proposition}
In fact, once we have Proposition \ref{pp:b2}, we can take a finite subset $\{y(j)\}_{j=1}^{J}$ in $S^1$ so that the neighborhoods $U_{y(j)}$ in Proposition \ref{pp:b2} cover $S^1$ and, letting  $\{\varphi_i\}_{j=1}^{m}$ be the union of $\{\varphi_{y(j),i}\}_{i=1}^{\ell^\nu}$ for $1\le j\le J$ in the corresponding conclusions of Proposition \ref{pp:b2}, we obtain Proposition \ref{pp:b}.

\begin{proof}[Proof of Proposition \ref{pp:b2}]
Take a point $y\in S^1$ arbitrarily. For $\a,\b\in \A^\nu$, we write 
$\a \prec \b$ if $\tau^q(\b(y))=\a(y)$ for some $q\ge 0$. By  simple combinatorial argument, we can show that this is a partial order on $\A^\nu$ and that, for each $\a\in \A^\nu$, there exists at most $(\nu+1)$ elements $\b\in \A^{\nu}$ such that $\b\prec \a$. (See the proof of \cite[Proposition 16]{T}.) 

For $0<\epsilon<1/2$ and $\a\in \A^\nu$, let $U(\epsilon)$ be the $\epsilon$-neighborhood of $y$ and $U_{\a}(\epsilon)$ the connected component 
of $\tau^{-\nu}(U(\epsilon))$ that contains $\a(y)$.
We consider an integer $\mu>\nu$ that will be specified later. 
We then choose $\epsilon_0>0$ so small that $\tau^i(U_{\b}(\epsilon_0))\cap U_{\a}(\epsilon_0) \neq \emptyset$ for some $1\le i\le \mu$ only if $\a\prec \b$. 
Take functions $\varphi_{\a}\in C^\infty(S^1)$ for $\a\in \A^\nu$ supported on 
$U_{\a}(\epsilon_0)$ such that 
\[
 \frac{d}{dx}\varphi_{\a}(y)=\ell^{\nu}\quad \mbox{ on $U_{\a}(\epsilon_0/3)$}\quad \mbox{ and  }\quad \left|\frac{d}{dx}\varphi_{\a}(y)\right|< 2\ell^{\nu}\quad \mbox{ on $S^1$.}
 \]
Finally let $\varphi_{y,i}$, $1\le i\le \ell^\nu$ be a rearrangement of $\varphi_{\a}$, $\a\in\A^\nu$  and let $U_y=U(\epsilon_0/3)$.  

We show that the conclusion of the proposition holds for the neighborhood $U_y$ and the functions $\varphi_{y,i}$, $1\le i\le \ell^\nu$, provided that the integer $\mu$ is sufficiently large. 
Consider the family (\ref{eqn:fam}) with $\varphi_i=\varphi_{y,i}$ and $m=\ell^\nu$ and suppose that a subset
$A=\{\a_i\}_{1\le i\le (\nu+1)(p+1)}$ of $\A^\nu$ is given. From the property of the partial order $\prec$ on $\A^\nu$ mentioned above, we can choose a subset $A'=\{\a'_i\}_{0\le i\le p}$ of $A$ that consists of maximal elements in $A$ with respect to $\prec$. Let $\sigma=\{\b_i\}_{0\le i\le p}$ be a subset of $\A^n$ with $n\ge \nu$ such that $[\b_i]_{\nu}=\a'_i$ for $0\le i\le p$.  For $\b\in \A^n$ and $x\in U_y$, we put 
\[
h_1(x,\b;\t)= \sum_{j=1}^{\min\{n,\mu\}}\ell^{-j} \frac{d}{dx}f_{\t}([\b]_{j}(x))
\]
and
\[
h_2(x,\b;\t)= \sum_{j=\min\{n,\mu\}+1}^{n}\ell^{-j} \frac{d}{dx}f_{\t}([\b]_{j}(x)),
\]
so that
\[
\frac{ds}{dx}(x, \b;f_\t)=h_1(x,\b;\t) +h_2(x,\b;\t).
\] 
Accordingly we decompose the affine map $G_{x,\sigma}$ into 
\[
G_{x,\sigma}^{(1)}(\t)=
\left(
h_1(x,\b_i;\t)-h_1(x,\b_0;\t)
\right)_{i=1,2,\dots,p}: \real^{\ell^\nu}\to \real^p
\]
and 
\[
G_{x,\sigma}^{(2)}(\t)=
\left(
h_2(x,\b_i;\t)-h_2(x,\b_0;\t)
\right)_{i=1,2,\dots,p}: \real^{\ell^\nu}\to \real^p.
\]
Let $\xi:\{1,2,\cdots,p\}\to \{1,2,\cdots,\ell^\nu\}$ 
be the correspondence such that  
$\a'_i(y)\in \supp(\varphi_{y,\xi(i)})$ for $1\le i\le p$, and consider the subspace of $\real^{\ell^\nu}$, 
\[
E=\{\t=(t_j)_{j=1}^{\ell^\nu}\in \real^{\ell^\nu}
\mid t_j\neq 0\mbox{ only if  $j=\xi(i)$ for some $1\le i\le p$}\},
\]
which is naturally identified with $\real^p$. 
Take any point $x\in U_y$ and let
$L^{(1)}$ and $L^{(2)}$ be the matrices that represent the linear part of the affine mappings 
$G_{x,\sigma}^{(1)}:E\to \real^p$ and $G_{x,\sigma}^{(2)}:E\to \real^p$ respectively. As a consequence of the choice of $\a'_i$, we can see that $L^{(1)}$ is the identity matrix of size $p$ while all the entries $L^{(2)}$ are bounded by $2\ell^{-\mu+\nu}(1-\ell^{-1})^{-1}$. 
Therefore, if we take sufficiently large $\mu$, it holds 
\[
\Jac(DG_{x,\sigma})\ge \Jac(DG_{x,\sigma}|_{E})\ge 1/2.
\]
Multiplying each $\varphi_{y,i}$ by $2$, we can replace $1/2$ by $1$ on the right hand side. 
\end{proof}

\appendix


\def\hcone{\hat{\cone}}

\section{Proof of Theorem \ref{th:tfae} }\label{sec:apd1}
We first show that the semi-flow $\mathbf{T}_f$ is weakly mixing if $\m(f)=1$. For this purpose, we introduce two quantities $\mathbf{n}(f,t)$ and $\mathbf{n}(f)$, similar to $\m(f,t)$ and $\m(f)$ respectively, as follows. 
Put $\hcone_f=\{(x,y)\in \real^2\mid |y|\le 2\theta_f|x|\}\supset \cone_f$. 
For $t\ge 0$, $z\in X_f$ and a one dimensional subspace $L\subset \real^2$, we define 
\[
\mathbf{n}(f,t,z,L)=
\sum{}^*\frac{1}{E(\zeta,t;f)}\le 1
\]
where $\sum^*$ is the sum over $\zeta\in (T^{t}_f)^{-1}(z)$ such that $(DT_f^t)_\zeta(\hcone_f)\supset L$. Then we put
\[
\mathbf{n}(f,t)=\max_{z\in X_f} \max_{L\in \real\mathbf{P}^1}\; \mathbf{n}(f,t, z, L)
\]
and 
\[
\mathbf{n}(f)=\limsup_{t\to \infty} \mathbf{n}(f,t)^{1/t}.
\]
Note that $\mathbf{n}(f,t)$ is {\em sub-multiplicative} with respect to $t$:
$\mathbf{n}(f,t+s)\le \mathbf{n}(f,t)\cdot \mathbf{n}(f,s)$. 
In this point, the quantity $\mathbf{n}(f,t)$ is better than $\m(f,t)$. In particular,  the limit in the definition of $\mathbf{n}(f)$ is actually exact.

We first show that $\m(f)=1$ implies $\mathbf{n}(f)=1$. 
For this purpose, it is sufficient to prove the claim that
\[
\m(f,s)\le \mathbf{n}(f,t)\quad \mbox{ for any $t\ge 0$ and $s=(b/a)t+b>t$}
\]
where  
\[
a=\min_{x\in S^1}f(x)\quad \mbox{and}\quad  b=\max_{x\in S^1}f(x).
\]
Consider a point $z\in X_f$ and take $w\in T^{-s}_f(z)$. 
If  
\begin{equation}\label{eqn:pit}
(DT^s_f)_\zeta(\cone_f)\cap (DT^s_f)_{w}(\cone_f)\neq \{0\}
\end{equation}
for a points $\zeta\in T^{-s}_f(z)$, then we have 
\begin{equation}\label{eqn:cl}
(DT^{t}_f)_{\zeta'}(\hcone_f)\supset L:=(DT^s_f)_w(\real\times \{0\})\quad \mbox{for $\zeta'=T^{s-t}_f(\zeta)\in T_f^{-t}(z)$.}
\end{equation}
Indeed, this follows from the fact that the differences between the slope of $L$ and those of boundary lines of $(DT_f^s)_{w}(\cone_f)$ are not greater than $\ell^{-[s/b]}\theta_f$, while the differences between the slopes of the boundary lines 
of $(DT_f^t)_{\zeta'}(\hcone_f)$ and those of the boundary lines 
of $(DT_f^t)_{\zeta'}(\cone_f)$ are  greater than $\ell^{-[t/a]-1}\theta_f= \ell^{-[s/b]}\theta_f$. Hence, in view of (\ref{sum}), we have that
\[
\sum_{\zeta:\zeta\pitchfork w} \frac{1}{E(\zeta,s;f)}
\le \sum{}^* \frac{1}{E(\zeta,t;f)}
\]
where $\sum_{\zeta:\zeta\pitchfork w}$ denotes the sum over $\zeta\in T^{-s}_f(z)$ satisfying (\ref{eqn:pit}) and $\sum^*$ denotes the same sum as that in the definition of $\mathbf{n}(f,t,z,L)$. Clearly this implies the claim above.

We next show that $\mathbf{T}_f$ is weakly mixing if $\mathbf{n}(f)=1$. Suppose  $\mathbf{n}(f)=1$. By submultiplicative property of $\mathbf{n}(f,t)$, we have  $\mathbf{n}(f,t)=1$ for all $t\ge 0$. Therefore we can take sequences of real numbers
 $t_n\ge 0$, points $z_n\in X_f$ and one-dimensional subspaces $L_n \in  \real\mathbf{P}^1$ for $n\ge 1$ such that
$\mathbf{n}(f,t_n, z_n, L_n)=1$ for all $n\ge 1$ and that, as $n\to \infty$,
\begin{itemize}
\item $t_n\to \infty$, 
\item $z_n$ converges to some $z_\infty\in X_f$,
and 
\item $L_n$ converges to some $L_\infty$ (in $\real\mathbf{P}^1$).
\end{itemize}
The condition $\mathbf{n}(f,t_n, z_n, L_n)=1$ and (\ref{sum})  imply that the cone $(DT^{-t_n}_f)_{w}(\hcone_f)$ contains $L_n$ for all the point $w\in T^{-t_n}_f(z_n)$. 
Note the unstable subspace (or the tangent space of the unstable manifold) for a backward orbit 
$(w(t))_{t\le 0}$ is contained in $(DT^{-t}_f)_{w(t)}(\cone_f)$ for any $t\le 0$. Thus, by continuity, we see that 
the unstable subspaces for {\em all} backward orbits of $z_\infty$ coincide with each other (and with $L_\infty$). 
Moreover, such property holds not only for the point $z_\infty$ but for all the points in $X_f$ because the set of points with such property is closed and completely invariant with respect to the flow $\mathbf{T}_f$.

For $x\in S^1$, let $\psi(x)$ be the slope of the (unique) ustable subspace at $(x,0)\in X_f$.
Invariance of the unstable subspaces implies that we have
 \begin{equation}\label{difcoc}
\psi(\tau(x))= (f'(x)+\psi(x))/\ell\quad \mbox{for all $x\in S^1$.}
\end{equation}
Inductive use of this equality yields
\begin{equation}\label{eqn:psie}
\psi(x)=\sum_{n\ge 1} \sum_{\tau^n(y)=x} \ell^{-2n}f'(y) \quad \mbox{for all $x\in S^1$,}
\end{equation}
where the right hand side converges in $C^{r-1}$ sense. 
Since we have $\int_{S^1} \psi(x)=0$ from (\ref{eqn:psie}), the function  
\[
\Psi(x)=\int_0^x \psi(y)dy
\]
is well-defined and $C^r$ on $S^1$. It follows from (\ref{difcoc}) that 
\[
\Psi(\tau(x))=\Psi(x) +f(x)-c\qquad \mbox{for some constant $c$.}
\]
By integrating the both sides over $S^1$, we see that $c=\int_{S^1} f(x) dx>0$. 
Now define   
\[
\Phi(x,s)=\exp((2\pi i/c)(\Psi(x)+ s))\quad \mbox{ for $(x,s)\in X_f$.}
\]
Then  $\Phi\circ T^t_f=e^{(2\pi i/c) t}\Phi$ for $t\ge 0$ . Therefore $T^t_f$ is not weakly mixing.

To finish the proof, we show that  $\m(f)=1$ if the semi-flow $\mathbf{T}^t_f$ is not weakly mixing. Suppose that $\mathbf{T}^t_f$ is not weakly mixing. Then we can find a real number $a\neq 0$ and an $L^2$ function $\Phi$ on $X_f$ such that $\Phi\circ T^t_f=e^{ia t}\Phi$ for $t\ge 0$. Equivalently there exists an $L^2$ function $\Psi$ on $S^1$ such that $\Phi(x,s)=e^{-ias} \Psi(x)$ and that $\Psi(\tau(x))=e^{iaf(x)}\Psi(x)$ for $x\in S^1$ and $t\ge 0$.
Actually the last equality tells that $\Psi$ is a $C^r$ function and so is $\Phi$. (For the proof of this fact, we refer that of \cite[Proposition 4.2]{PP}\footnote{Replace the symbolic dynamical system $\sigma:X^+\to X^+$ and the space of H\"older functions on $X^+$ in the proof of \cite[Proposition 4.2]{PP}  by  $\tau:S^1\to S^1$ and $C^r(S^1)$ respectively.}, for instance.) Let $L(z)$ be the null line of the differential $D_z\Phi$. Then this line field is invariant with respect to the semi-flow $\mathbf{T}_f$ and not tangent to the flow direction. Hence $L(z)$ is contained in the cone $(DT^{-t}_f)_{w(t)}(\cone_f)$ for any backward orbit $\{w(t)\}_{t\le 0}$ of $z$ and any $t\le 0$. This and (\ref{sum}) imply that $\m(f,t)=1$ for any $t\ge 0$ and hence that $\m(f)=1$.

\section{Proof of Lemma \ref{lm:LYS}}\label{sec:apd2} 
Let $\Gamma=\mathbb{Z}_+\times \{+,-\}$, $c(+)=1$ and $c(-)=0$.
Below we write $C_0$ for constants that does not depend on $S$, $h$, $\Theta$ nor $\Theta'$, while we write $C$ for constants that may depend on them.  Take an integer $\mu=\mu(S)$ such that 
\[
2^{-\mu+6}\|\xi\|\le \|(DS_{x})^{tr}(\xi)\|\le 2^{\mu-6}\|\xi\|
\qquad\mbox{for  any  $x\in K$ and any $\xi\in \real^2$.}
\]
Let $\nu\le \mu-6$ be an integer such that 
\[
2^{\nu-1}< \Lambda(S,\Theta',K)\le 2^{\nu}.
\] 
So we have 
\[
\|DS^{tr}_{x}(\xi)\|\le 2^{\nu}\|\xi\|
\qquad \mbox{ if $x\in K$ and $(DS_{x})^{tr}(\xi)\notin \cone'_{-}$.}
\] 
We write $(m, \tau) \hookrightarrow (n,\sigma)$ if either 
\begin{itemize}
\item $(\tau,\sigma)=(+,+)$ and $m-\mu\le n\le \max\{0,m+\nu+6\}$, or 
\item $(\tau,\sigma)\in \{(-,-), (+,-)\}$  and 
$m-\mu\le n\le m+\mu$.
\end{itemize}
And we write $(m, \tau) \not\hookrightarrow (n,\sigma)$ otherwise. 

Consider a function $u\in C^r(R)$ and put  $v:=L u$. For $(n,\sigma), (m, \tau)\in \Gamma$, we define 
\[
v_{n,\sigma}^{m,\tau}=\psi_{\Theta',n,\sigma}(D) 
L (u_{\Theta, m, \tau}),
\] so that $v_{\Theta',n,\sigma}=\sum_{(m,\tau)\in \Gamma} v_{n,\sigma}^{m,\tau}$. 
By using Parseval's identity, we can get
\begin{equation}\label{eqn:simple}
\sum_{(n,\sigma)\in \Gamma} 
\|v_{n,\sigma}^{m,\tau}\|_{L^2}^2 \le \left\|L (u_{\Theta, m, \tau})\right\|_{L^2}^2
\le C_0 \gamma(S)^{-1}
\|h\|_{L^{\infty}}^2\|u_{\Theta,m,\tau}\|_{L^2}^2.
\end{equation}
Also we have the following estimate, whose proof is postponed for a while.  
\begin{lemma}\label{lm:s}
If $(m,\tau)\not\hookrightarrow(n, \sigma)$, we have 
\begin{equation}\label{eqn:s}
\|v_{n,\sigma}^{m,\tau}\|_{L^2}\le C
2^{-(r-1)\max\{m,n\}}\|u_{\Theta, m, \tau}\|_{L^2}.
\end{equation}
\end{lemma}
\begin{remark}
If $S$ is an affine map in the lemma above, the Fourier transform of $L (u_{\Theta, m, \tau})$ is supported on 
$DS^{tr}(\supp(\psi_{\Theta,m,\tau}))$, which does not meet $\supp(\psi_{\Theta',n,\sigma})$ by the assumption $(m,\tau)\not\hookrightarrow(n, \sigma)$,  and hence the assertion holds trivially with $v_{n,\sigma}^{m,\tau}=0$.
To prove the lemma above, we will estimate some oscillatory integrals using smoothness of $S$. 
\end{remark}
We first show the assertion that $|v|_{\Theta'}\le C|u|_{\Theta}$ for some constant $C$, which implies that $L$ extends boundedly to $L:W_\dag(R;\Theta)\to W_\dag(R;\Theta')$.  
By definition, we have 
\[
|v|_{\Theta'}^2\le 2
\sum_{(n,\sigma)\in \Gamma}2^{2(c(\sigma)-\epsilon)n}
\left(
\left\|
\sum_{(m,\tau)\hookrightarrow (n,\sigma)}  v_{n,\sigma}^{m,\tau}
\right\|_{L^2}^2
+\left\|
\sum_{(m,\tau)\not\hookrightarrow (n,\sigma)} v_{n,\sigma}^{m,\tau}
\right\|_{L^2}^2\right)
\]
where $\sum_{(m,\tau)\hookrightarrow (n,\sigma)}$ (resp. $\sum_{(m,\tau)\not\hookrightarrow (n,\sigma)}$) denotes the sum over $(m,\tau)\in \Gamma$ such that $(m,\tau)\hookrightarrow (n,\sigma)$ (resp. $(m,\tau)\not\hookrightarrow (n,\sigma)$). 
Since the relation $(m,\tau)\hookrightarrow (n,\sigma)$ holds only if $c(\sigma)\le c(\tau)$ and $|m-n|<\mu$, it holds
\begin{align*}
\sum_{(n,\sigma)\in \Gamma}
\left\|
\sum_{(m,\tau)\hookrightarrow (n,\sigma)} 2^{(c(\sigma)-\epsilon)n} v_{n,\sigma}^{m,\tau}
\right\|_{L^2}^2& \le 
C\sum_{(n,\sigma)\in \Gamma}
\sum_{(m,\tau)\in \Gamma} 
2^{2(c(\tau)-\epsilon)m}\| v_{n,\sigma}^{m,\tau}\|_{L^2}^2\\
&\le C
\sum_{(m,\tau)\in \Gamma} 
2^{2(c(\tau)-\epsilon)m}\|u_{\Theta,m,\tau}\|_{L^2}^2\le C
|u|_{\Theta}^2
\end{align*}
where the second inequality follows from (\ref{eqn:simple}). 
Also it follows from Lemma \ref{lm:s} and Schwarz inequality  that
\begin{align}
&\sum_{(n,\sigma)\in \Gamma}
\left\|
\sum_{(m,\tau)\not\hookrightarrow (n,\sigma)} 2^{(c(\sigma)-\epsilon)n} v_{n,\sigma}^{m,\tau}
\right\|_{L^2}^2\label{eqn:uTheta}\\
&\le \sum_{(n,\sigma)\in \Gamma}
\left\|
\sum_{(m,\tau)\not\hookrightarrow (n,\sigma)} \!\!\!\!\!\!\!
2^{c(\sigma)n-c(\tau)m-(r-1-\epsilon)\max\{n,m\}}\cdot 2^{(r-1)\max\{n,m\}+(c(\tau)-\epsilon)m} v_{n,\sigma}^{m,\tau}
\right\|_{L^2}^2 \notag\\
&\le 
\sum_{(n,\sigma)\in \Gamma}\!\!
\left(
\sum_{(m,\tau)\in\Gamma}\!\!\!\!
2^{2c(\sigma)n-2c(\tau)m-2(r-1-\epsilon)\max\{n,m\}}
\right)\!\!
\left(
\sum_{(m,\tau)\in \Gamma}\!\!\!\!
2^{2(c(\tau)-\epsilon)m}\|u_{\Theta,m,\tau}\|_{L^2}^2
\right)\notag\\
&\le C|u|_{\Theta}^2.\notag
\end{align}
Thus we obtain $|v|_{\Theta'}\le C |u|_{\Theta}$ for $u\in C^r(R)$ and hence for $u\in W_\dag(R;\Theta)$. 

We next prove (\ref{eqn:LYSm}) and (\ref{eqn:LYSp}). 
The inequality (\ref{eqn:LYSm}) is easy to see: 
\begin{align*}
(\|v\|_{\Theta'}^-)^2
\le  \|v\|_{L^2}^2
\le \gamma(S)^{-1} \|h\|_{L^{\infty}}^2\|u\|_{L^2}^2
\le \gamma(S)^{-1} \|h\|_{L^{\infty}}^2\|u\|_{\Theta}^2. 
\end{align*} 
To prove (\ref{eqn:LYSp}), we begin with writing the left hand side as 
\[
(\|v\|_{\Theta'}^+)^2=\|\psi_{\Theta',0,+}(D)v\|_{L^2}^2
+\sum_{n\ge 1}2^{2n}\|\psi_{\Theta',n,+}(D)v\|_{L^2}^2.
\]
The first term on the right hand side is bounded by 
$|v|_{\Theta'}^2$ and hence by $C|u|^2_{\Theta}$.
The sum on the right hand side 
is bounded by  
\[
2\cdot \left( \sum_{n\ge 1}
\left\|
\sum_{(m,\tau)\hookrightarrow (n,+)} 2^{n} v_{n,+}^{m,\tau}
\right\|_{L^2}^2
+
\sum_{n\ge 1}
\left\|
\sum_{(m,\tau)\not\hookrightarrow (n,+)} 2^{n} v_{n,+}^{m,\tau}
\right\|_{L^2}^2\right).
\]
By Schwarz inequality, we have
\begin{align*}
\left\|
\sum_{(m,+)\hookrightarrow (n,+)} 2^{n} v_{n,+}^{m,+}
\right\|_{L^2}^2
 \le  \left(\sum_{(m,+)\hookrightarrow (n,+)} 2^{2(n-m)} \right)
\left(\sum_{(m,+)\hookrightarrow (n,+)}
2^{2m}\|v_{n,+}^{m,+}\|_{L^2}^2\right)
\end{align*}
where $\sum_{(m,+)\hookrightarrow (n,+)}$ denotes the sum over $m\ge 0$ such that $(m,+)\hookrightarrow (n,+)$. 
Note that we have $(m,\tau)\hookrightarrow (n,+)$ for $n\ge 1$ only if  $\tau=+$ and $n\le m+\nu+6$. Thus we can see, by using (\ref{eqn:simple}), that 
\begin{align*}
\sum_{n\ge 1}\left\|
\sum_{(m,\tau)\hookrightarrow (n,+)} 2^{n} v_{n,+}^{m,\tau}
\right\|_{L^2}^2
&\le C_0 \cdot 2^{2\nu}\gamma(S)^{-1} \|h\|_{L^{\infty}}^2\|u\|_{\Theta}^2.
\end{align*}
On the other hand, by using Lemma \ref{lm:s} and Schwarz inequality, we can show
\[
\sum_{n\ge 0}\left\|
\sum_{(m,\tau)\not\hookrightarrow (n,+)} 2^{n} v_{n,+}^{m,\tau}
\right\|_{L^2}^2
<C|u|_{\Theta}^2
\]
in the similar manner as (\ref{eqn:uTheta}).
Therefore we obtain (\ref{eqn:LYSp}). Obviously (\ref{eqn:LYSm}) and (\ref{eqn:LYSp}) imply (\ref{eqn:LYSpm}) and hence $L$ extends boundedly to $L:W_*(R;\Theta)\to W_*(R;\Theta')$. 
Finally we complete the proof by proving Lemma \ref{lm:s}.

\begin{proof}[Poof of Lemma \ref{lm:s}]
Since $K$ is compact, we can take closed cones $
\widetilde \cone_+ \Subset \cone_+$ and $
\;\widetilde \cone_-\Subset \cone_-$ such that 
\[
(DS_\zeta)^{tr}(\real^d\setminus \widetilde \cone_{+}) \Subset  \cone'_{-}
\quad \mbox{ for $\zeta\in K$. }
\]
Let 
$\tilde \varphi_+$, $\tilde\varphi_-:S^1\to [0,1]$ be
$C^{\infty}$ functions  satisfying
\[
\tilde \varphi_+(\xi)=
\begin{cases}
1, &\mbox{if $\xi\notin S^1 \cap \cone_{-}$;}\\
0, &\mbox{if $\xi\in S^1 \cap \widetilde \cone_{-}$,}
\end{cases} \qquad 
\tilde\varphi_-(\xi)=
\begin{cases}
1, &\mbox{if $\xi\notin S^1 \cap \cone_{+}$;}\\
0, &\mbox{if $\xi\in S^1 \cap \widetilde\cone_{+}$.}
\end{cases}
\]
Recall the function $\chi$ and define $
\tilde\psi_n(\xi)=
\chi(2^{-n-1}|\xi|)-\chi(2^{-n+2}|\xi|)$ for $n\ge 1$ and
\[
\tilde \psi_{\Theta, n, \sigma}(\xi)=
\begin{cases}
\tilde{\psi}_n(\xi)\tilde \varphi_{\sigma}(\xi/|\xi|),& \mbox{ if $n\ge 1$;}\\
\chi(2^{-1}|\xi|),& \mbox{ if $n=0$}
\end{cases}
\]
for $(n,\sigma)\in\Gamma$.  
Then we have $\tilde \psi_{\Theta,n,\sigma}(\xi)=1$ 
if $\xi \in \supp (\psi_{\Theta,n,\sigma})$. From the definition of the relation $\hookrightarrow$,   
there exists a constant $L>1$, which may depend on $S$,  such that, if $(m, \tau) \not\hookrightarrow (n,\sigma)$ and $\max\{m,n\}\ge L$, it holds
\begin{equation}\label{lowerbd}
d(\supp(\psi_{\Theta',n,\sigma}), (DS_\zeta)^{tr}(\supp(\tilde \psi_{\Theta,m,\tau})))
\ge L^{-1}\cdot 2^{\max\{n,m\}} 
\quad \mbox{for  $\zeta\in K$.}
\end{equation}
In the case where $\max\{m,n\}< L$, it is easy to see that (\ref{eqn:s}) holds with the constant $C$ depending on $L$. 
Thus we assume $\max\{m,n\}\ge L$ in the following.

We consider the operator $S_{n,\sigma}^{m,\tau}$ defined by 
\[
S_{n,\sigma}^{m,\tau}= \psi_{\Theta',n,\sigma}(D)\circ L\circ \tilde{\psi}_{\Theta,m,\tau}(D).
\]  
Then we have $v_{n,\sigma}^{m,\tau}=S_{n,\sigma}^{m,\tau} u_{\Theta,m,\tau}$ because $\tilde{\psi}_{\Theta,m,\tau}(D)(u_{\Theta,m,\tau})=u_{\Theta,m,\tau}$. 
We may rewrite this operator $S_{n,\sigma}^{m,\tau}$ as 
\[
(S_{n,\sigma}^{m,\tau}u)(x)
=(2\pi)^{-4}\int V_{n,\sigma}^{m,\tau}(x,y) \cdot u\circ S(y) \cdot |\det DS(y)| 
dy,
\]
where 
\begin{equation}\label{Vkernel}
V_{n,\sigma}^{m,\tau}(x,y)=\int e^{i(x-w)\xi+i(S(w)-S(y))\eta} 
h(w)\psi_{\Theta',n,\sigma}(\xi)
\tilde{\psi}_{\Theta,m,\tau}(\eta)dw d\xi d\eta.
\end{equation}
Since we have $\| u\circ S(y) \cdot |\det DS(y)| \|_{L^2}\le C\|u\|_{L^2}$, 
the inequality (\ref{eqn:s}) 
follows if 
the operator norm of the integral operator 
\[
H^{m,\tau}_{n,\sigma}:L^2(\real^2)\to L^2(\real^2),  
\quad H^{m,\tau}_{n,\sigma}v(x)= \int V_{n,\sigma}^{m,\tau}(x,y)v(y) dy
\]
is 
bounded by $C\cdot 2^{-(r-1)\max\{n,m\}}$.

Apply the following formula of integration by parts  
for $(r-1)$ times in (\ref{Vkernel}),
\[
 \int e^{if(w)}g(w) dw=
i\cdot \int e^{if(w)}\cdot \sum_{k=1}^2 \partial_{w_k}\left(\frac
{\partial_{w_k} f(w)\cdot  g(w)}
{\sum_{j=1}^{2}(\partial_{w_j} f(w))^2}\right) dw
\]
where $w=(w_k)_{k=1}^{2}\in \real^2$. Then we obtain the expression 
\[
V_{n,\sigma}^{m,\tau}(x,y)=
\int e^{i(x-w)\xi+i(S(w)-S(y))\eta} 
F(\xi,\eta,w)\psi_{\Theta',n,\sigma}(\xi)
\tilde{\psi}_{\Theta,m,\tau}(\eta)dwd\xi d\eta
\]
where $F(\xi,\eta,w)$ is continuous in $w$ and $C^\infty$ in $\xi$ and $\eta$. Note that
$F(\xi,\eta,w)=0$ if $w\notin K$. From (\ref{lowerbd}), there is a constant $C_{\alpha\beta}$ for multi-indices $\alpha$ and $\beta$, such that 
\begin{equation}\label{eqn:scf}
|\partial_\xi^\alpha\partial_\eta^\beta
F(\xi,\eta,w)|\le C_{\alpha\beta}\cdot 
2^{-n|\alpha|-m|\beta|-(r-1)\max\{n,m\}}
\end{equation}
for $w\in \real^2$, $\xi\in \supp(\psi_{\Theta',n,\sigma})$ and $\eta\in \supp(\tilde{\psi}_{\Theta,m,\tau})$.
For $n\ge 0$ and $m\ge 0$, we put 
\[
G_{nm}(\xi,\eta,w)= F(2^{n}\xi,2^{m}\eta,w)\psi_{\Theta',n,\sigma}(2^{n}\xi)
\tilde{\psi}_{\Theta,m,\tau}(2^{m}\eta).
\]
By changes of variable,  we can rewrite $V_{n,\sigma}^{m,\tau}(x,y)$ as
\begin{equation}\label{eqn:v}
\int 2^{2n+2m}
(\mathcal{F}_{\xi\eta}^{-1}G_{nm})(2^{n}(x-w),2^{m}(S(w)-S(y)),w) dw
\end{equation}
where $\mathcal{F}_{\xi\eta}^{-1}$ is the inverse Fourier transform with respect to the variables $\xi$ and $\eta$. 
From (\ref{eqn:scf}), there exists a constant $C_{\alpha\beta}$ for any multi-indices $\alpha$ and $\beta$ such that 
\[
|\partial_\xi^\alpha\partial_\eta^\beta G_{nm}|_{L^\infty}\le C_{\alpha\beta} 2^{-(r-1)\max\{n,m\}}.
\]
This implies that there exists a constant $C$ such that
\[
|\mathcal{F}_{\xi\eta}^{-1}G_{nm}(x,y,w)|\le C 2^{-(r-1)\max\{n,m\}}(1+|x|^2)^{-2}(1+|y|^2)^{-2}.
\]
Applying this inequality in the expression (\ref{eqn:v}) of $V_{n,\sigma}^{m,\tau}(x,y)$, we obtain
 the required estimate for  $H_{n,\sigma}^{m,\tau}$ from Young's inequality. 
\end{proof}

\section{Proof of Lemma \ref{lm:pu}}\label{sec:apd3}
We prove the inequality (\ref{eqn:pu1}). For the inequality (\ref{eqn:pu2}), we refer \cite[Lemma 7.1]{BT}. 
Recall the argument in the proof of Lemma~\ref{lm:LYS} in Appendix \ref{sec:apd2}, setting $S=id$. Notice that the assumptions of Lemma~\ref{lm:LYS} then hold since we assume $\Theta'<\Theta$. 
Put $v_i=g_i\cdot u$ for $1\le i\le I$. Then we have
\begin{equation}\label{eqn:trb}
\sum_{i=1}^I (\|v_i\|_{\Theta'}^-)^2\le
\sum_{i=1}^I \|v_i\|_{L^2} ^2\le 
\|u\|_{L^2}^2\le C \|u\|_{\Theta}^2.
\end{equation}
We have proved, in the proof of Lemma~\ref{lm:LYS},  that 
\[
\sum_{n\ge 0}2^{2n}
\left\|\sum_{(m,\tau)\not\hookrightarrow (n,+)}
\psi_{\Theta',n,+}(D)(g_i u_{\Theta, m,\tau})\right\|_{L^2}^2
\le 
C|u|_{\Theta}^2
\]
Since we can and do put $\mu=6$ in the setting $S=id$, the relation $(m,\tau)\hookrightarrow (n,+)$ holds only if $|m-n|\le 6$. Hence we have, by Schwarz inequality,
\begin{align*}
&\sum_{i=1}^{I}\sum_{n\ge 0}2^{2n}
\left\|\sum_{(m,\tau)\hookrightarrow (n,+)}
\psi_{\Theta',n,+}(D)(g_i u_{\Theta, m,\tau})\right\|_{L^2}^2\\
&\qquad\le 
13\cdot \sum_{i=1}^I\sum_{n\ge 0}\;\;\sum_{m:|m-n|\le 6}2^{2n}
\left\|
\psi_{\Theta',n,+}(D)(g_i u_{\Theta, m,+})\right\|_{L^2}^2 \\
&\qquad \le 13\cdot 2^{12}\cdot  \sum_{m\ge 0}\sum_{i}2^{2m}
\left\|g_i u_{\Theta, m,+} \right\|_{L^2}^2 \le C_0 \|u\|_{\Theta}^2.
\end{align*}
Therefore we obtain $
\sum_{i=1}^I (\|v_i\|_{\Theta'}^+)^2\le C_0 \|u\|_{\Theta}^2+C|u|_{\Theta}^2$. This and (\ref{eqn:trb}) yield (\ref{eqn:pu1}). 
\bibliographystyle{amsplain}

\begin{thebibliography}{10}
\bibitem{AGT}
A.~Avila, S.~Gou\"ezel \& M.~Tsujii,
\textit{Smoothness of solenoidal attractors,}
Discrete and Continuous Dynamical Systems {\bf 15} (2006), no. 1, 21--35.
\bibitem{VB} V.~Baladi, Positive transfer operators and decay of correlations. Advanced Series in Nonlinear Dynamics, 16. World Scientific, 2000
\bibitem{BT} V.~Baladi \& M.~Tsujii,
\textit 
{Anisotropic H\"older and Sobolev spaces 
for hyperbolic diffeomorphisms,} Ann. Inst. Fourier
 {\bf 57}(2007), no. 1, 127--154 
\bibitem{BV} V. Baladi \& B. Vall\'ee, \textit{Exponential decay of correlations for surface semi-flows without finite Markov partitions}. Proc. Amer. Math. Soc. {\bf 133} (2005), no. 3, 865--874
\bibitem{B} R. Bowen, Equilibrium states and the ergodic theory of Anosov diffeomorphisms, Lecture note in Math. {\bf 470}, Springer, 1975
\bibitem{Bo} R. Bowen, \textit{Periodic orbits for hyperbolic flows}. American J. Math. {\bf  94} (1972), 1-30

\bibitem{Do} D.~Dolgopyat, \textit{On decay of correlations in Anosov flows}, Ann.~Math.~{\bf 147} (1998), 357--390,
\bibitem{Do1} D.~Dolgopyat, \textit{Prevalence of rapid mixing in hyperbolic flows}, Ergodic Theory Dynam. Systems {\bf 18} (1998), no. 5, 1097--1114.
\bibitem{Do2} D.~Dolgopyat, \textit{Prevalence of rapid mixing in hyperbolic flows II, topological prevalence.} Ergodic Theory Dynam. Systems {\bf 20} (2000), no. 4, 1045--1059.
\bibitem{GL1} S. Gou\"ezel and C. Liverani,  \textit{Banach spaces adapted to Anosov systems}, Ergodic Theory Dynami. Systems, {\bf 26}(2006),  no. 1, 189--217 
\bibitem{GL2} S. Gou\"ezel and C. Liverani, Compact locally maximal hyperbolic sets for smooth maps: fine statistical properties, preprint(2006)
\bibitem{HSY}
B.R.\ Hunt, T.\ Sauer \& J.A.\ Yorke, \textit{Prevalence: a translation-invariant "almost every" on infinite-dimensional spaces}, Bull.\ Amer.\ Math.\ Soc.,  {\bf 27}, no.\ 2, 217--238(1992); Addendum {\it ibid.} {\bf 28}, no.\ 2, 306-307, (1993)

\bibitem{He} H. Hennion,
\textit{ Sur un th\'eor\`eme spectral
et son application aux noyaux lipschitziens,}
Proc. Amer. Math. Soc.  {\bf 118} (1993), 627--634.
\bibitem{IM} C. T. Ionescu Tulcea \& G.  Marinescu, 
\textit{Th\'eorie ergodique pour des classes d'op\'erations non compl\`etement continues. }
Ann. of Math. (2) {\bf 52} (1950), 140--147. 
\bibitem{L} C. Liverani, \textit{On contact Anosov flows}, Ann.~Math.~{\bf 159} (2004), 1275--1312.
\bibitem{PP} W. Parry \& M. Pollicott, \textit{Zeta functions and the periodic orbit structure of hyperbolic dynamics}, Ast\'erisque {\bf 187-188} (1990)
\bibitem{Po} M. Pollicott, \textit{On the mixing of Axiom A attracting flows and a conjecture of Ruelle}.  
Ergodic Theory Dynam. Systems {\bf 19} (1999), no. 2, 535--548.
\bibitem{R} D. Ruelle, A measure associated with Axiom A attractors, Amer. J. Math. {\bf 98} (1976), 616--654
\bibitem{S} Ya. G. Sinai, \textit{Gibbs measures in ergodic theory}, Russ. Math. Surveys {\bf 27} (1972), 21--70
\bibitem{T} M. Tsujii, \textit{Fat solenoidal attractors},  Nonlinearity {\bf 14} (2001), no. 5, 1011--1027.
\bibitem{T2}
M. Tsujii, \textit{A measure on the space of smooth mappings and dynamical system theory}, J.\ Math.\ Soc.\ Japan {\bf 44} (1992), no. 3, 415--425
\end{thebibliography}

\end{document}